\documentstyle[amstex,epsfig,amssymb,float,afterpage,silence]{amsart}
 \numberwithin{equation}{section}
 \newcommand{\abs}[1]{\left\vert#1\right\vert}

\newcommand{\br}[1]{\left( #1 \right)}
\newcommand{\Q}{\mathbb{Q}}
\newcommand{\R}{\mathbb{R}}

\newcommand{\F}{\mathbb{F}}

\begin{document}

\title[High Rank]
{Elliptic curves of high rank and the Riemann zeta function on the one line}

\author[Michael Rubinstein]
{M.\ O.\ Rubinstein\\ \\
\ Pure Mathematics \\ University of Waterloo\\ 200 University Ave W\\Waterloo, ON, Canada\\N2L 3G1\\}

%\subjclass{Primary L-functions}
%\keywords{L-functions, rank}

%\date{2008}

%%% ----------------------------------------------------------------------

\begin{abstract}
We describe some experiments that show a connection between elliptic curves
of high rank and the Riemann zeta function on the one line. We also
discuss a couple of statistics involving $L$-functions where the zeta function
on the one line plays a prominent role.
\end{abstract}

%%% ----------------------------------------------------------------------
\maketitle
%%% ----------------------------------------------------------------------

\section{Elliptic curves of high rank}

Let $E$ be an elliptic curve over $\Q$, which we write in minimal Weierstrass form
$$
   E: y^2 + a_1xy + a_3y = x^3 +a_2x^2 +a_4x +a_6,
$$
and denote by $E=[a_1,a_2,a_3,a_4,a_6]$.

To the elliptic curve $E$ we may associate an Euler product
\begin{eqnarray}
   \label{eq:euler_prdct}
    L_{E}(s)
    &=&\prod_{p\mid N}
    \left(1-a(p) p^{-(s+1/2)}\right)^{-1}
    \prod_{p\nmid N}
    \left(1-a(p) p^{-(s+1/2)}+p^{-2s}\right)^{-1}
    \notag \\
    &=&
    \prod_{p} {\mathcal{L}}_p(1/p^s)
    , \quad \quad \Re(s) > 1,
\end{eqnarray}
where $N$ is the conductor of $E$. For $p \nmid N$, $a(p) = p+1-\#E_p({\F}_p)$,
with $\#E_p({\F}_p)$ being the number of points on $E$ over ${\F}_p$. When
$p|N$, $a(p)$ is either $1$, $-1$, or $0$.

A theorem of Hasse states that $\abs{a(p)} < 2 p^{1/2}$.
Hence,~(\ref{eq:euler_prdct}) converges when $\Re(s)>1$, and for
these values of $s$ we may expand $L_E(s)$ in an absolutely convergent
Dirichlet series
\begin{equation}
    \label{eq:dirichlet series E}
    L_E(s)=\sum_1^\infty \frac{a(n)}{n^{1/2}} \frac{1}{n^s}.
\end{equation}

$L_E(s)$ has analytic continuation to $\mathbb{C}$ and satisfies a
functional equation of the form
\begin{equation}
    \label{eq:L_E functional eqn}
    \left(\frac{\sqrt{N}}{2\pi}\right)^{s}
    \Gamma(s+1/2)L_E(s)=w_E\left(\frac{\sqrt{N}}{2\pi}\right)^{1-s}\Gamma(3/2-s)L_E(1-s),
\end{equation}
where $w_E=\pm 1$ \cite{W} \cite{TW} \cite{BCDT}.

%Hence  we have
%\begin{align}
%    \br{\frac{N^{1/2}}{2\pi}}^{s} \Gamma(s+1/2) L_E(s)
%    \delta^{-s} =
%    \br{\frac{2 \pi \delta}{N^{1/2}}}^{1/2}
%    &\sum_{n=1}^{\infty}
%    a_n G\br{s+1/2, 2 \pi n \delta / N^{1/2}}
%    \notag \\
%    -\frac{\varepsilon}{\delta}
%    \br{\frac{2\pi}{N^{1/2} \delta}}^{1/2}
%    &\sum_{n=1}^{\infty}
%    a_n G\br{1-s+1/2, 2 \pi n /(\delta N^{1/2})}.
%    \notag
%\end{align}

We will describe some computational experiments involving
the elliptic curves with smallest known conductor of ranks 1--7, and rank 11~\cite{EW}:
\begin{table}[H]
\centerline{
\begin{tabular}{|c|c|c|}
\hline
$r$ (rank) & $E_r$ & $N$\cr
\hline
1 &$[0, 0, 1, -1, 0]$ & 37 \cr
2 &$[0, 1, 1, -2, 0]$ & 389 \cr
3 &$[0, 0, 1, -7, 6]$ & 5077 \cr
4 &$[1, -1, 0, -79, 289]$ & 234\,446  \cr
5 & $[0, 0, 1, -79, 342]$ & 19\,047\,851  \cr
6 & $[1, 1, 0, -2582, 48720]$ & 5\,187\,563\,742  \cr
7 & $[0, 0, 0, -10012, 346900]$ &382\,623\,908\,456  \cr
11 & $[0, 0, 1, -16359067, 26274178986]$ & 18\,031\,737\,725\,\\
& & $\qquad 935\,636\,520\,843$ \cr
\hline
\end{tabular}
}
\end{table}
as well
as an elliptic curve of rank at least 24 of Martin and McMillen~\cite{MM},
\begin{eqnarray}
    E_{24}=[1, 0, 1, -120039822036992245303534619191166796374, \notag \\
    504224992484910670010801799168082726759443756222911415116].
\end{eqnarray}

The above notation, of the form $E_r$, uses the rank of the elliptic curve as the
identifying subscript for the sake of compactness, rather than the usual way of
notating the elliptic curve in terms of its conductor (and more detailed information concerning its
isogeny class, which is not available here for $r\geq 5$).

\begin{figure}[H]
    \centerline{
            \psfig{figure=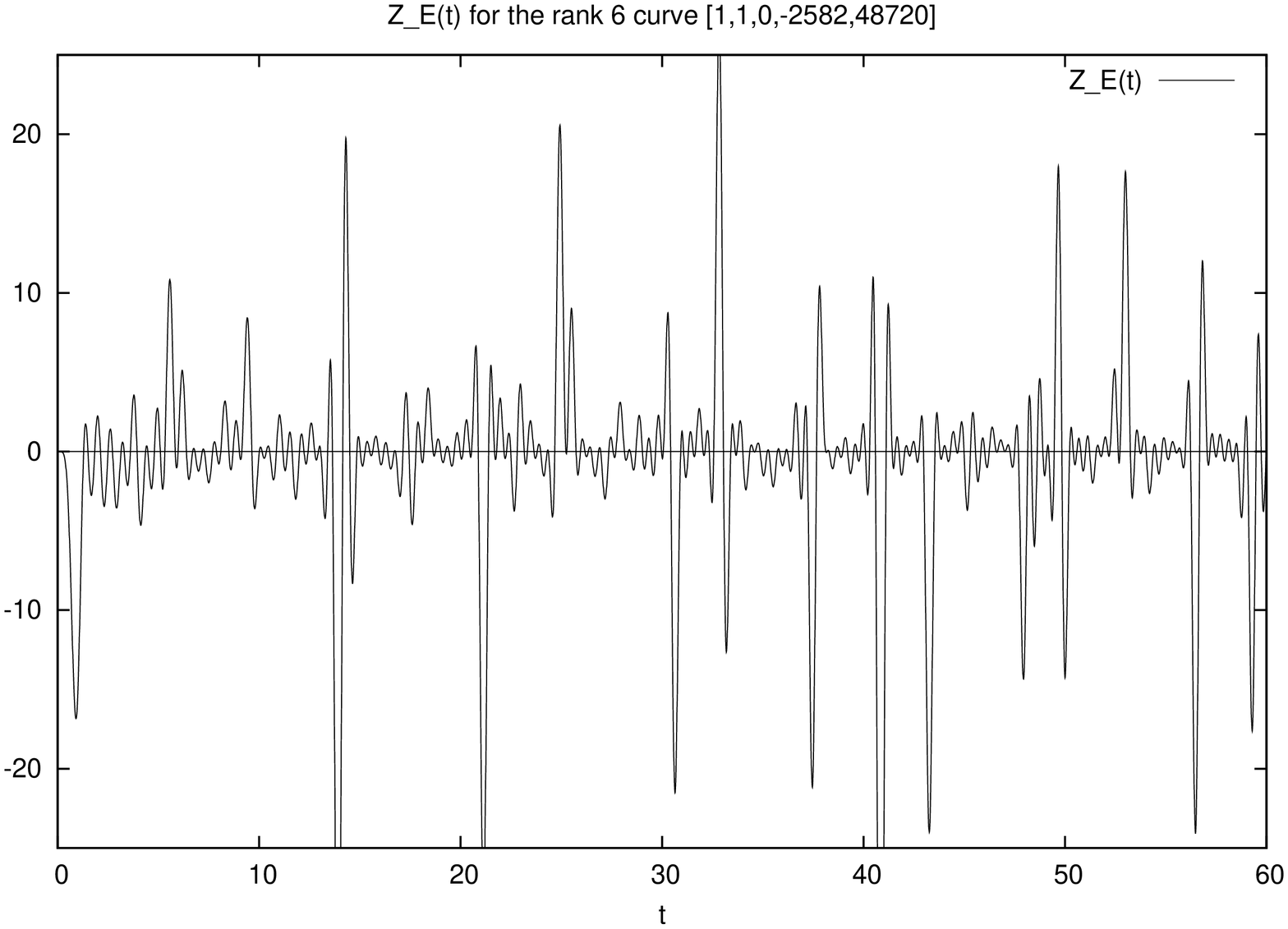,width=3in,angle=-90}
    }
    \centerline{
            \psfig{figure=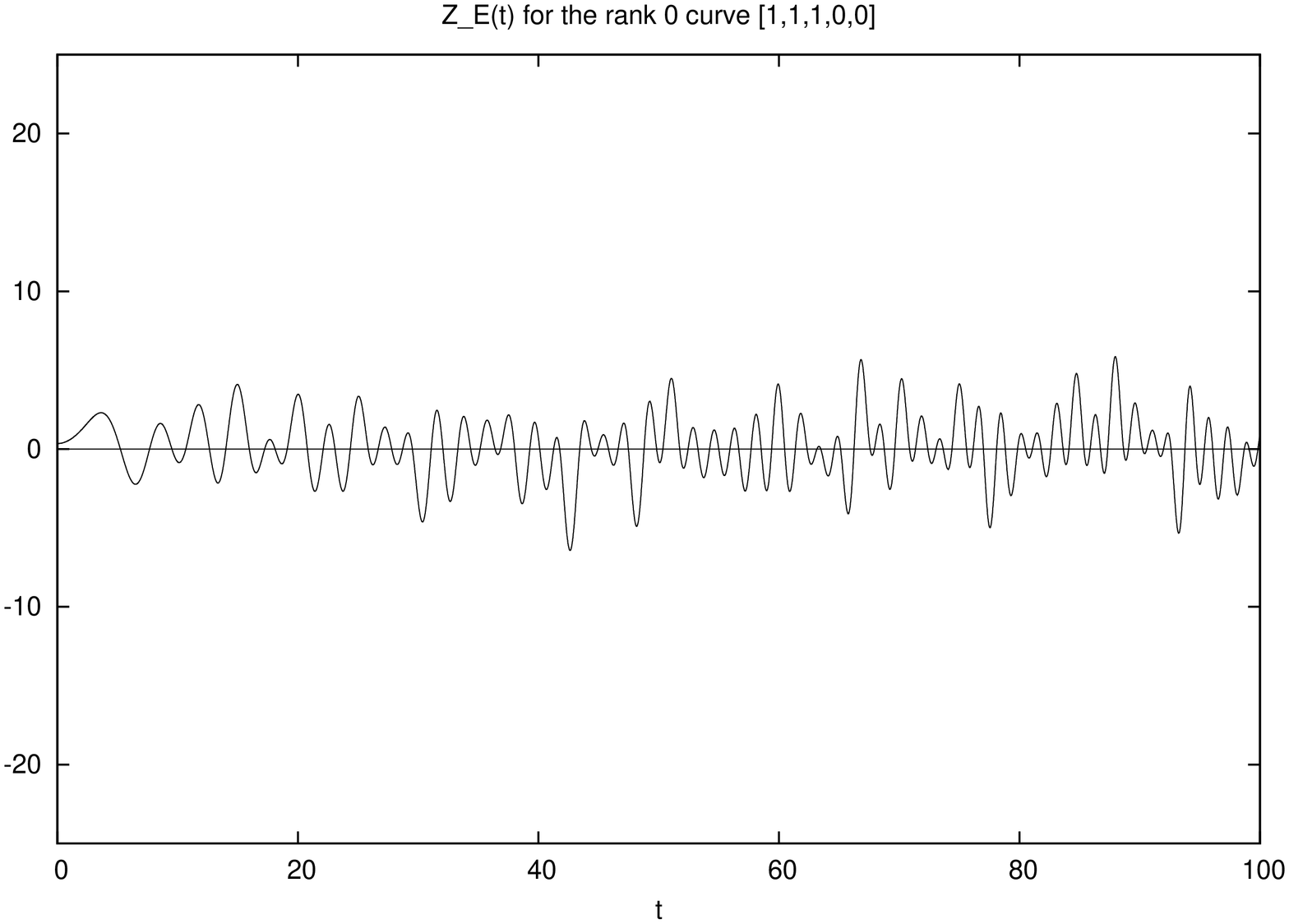,width=3in,angle=-90}
    }
    \caption
    {
        A plot (top) of the Hardy-$Z$ function of the $L$-function associated
        to the first known elliptic curve $[1,1,0,-2582,48720]$ of rank 6 and
        conductor $5\,187\,563\,742$, and (bottom) for the elliptic curve
        $[1,1,1,0,0]$ of rank 0 and conductor 15, $[1,1,1,0,0]$.
    }
    \label{fig:rank 6}
\end{figure}

While carrying out a numerical verification of the Riemann Hypothesis
for several moderately high rank elliptic curves, the author noticed peculiar behaviour when
plotting these $L$-functions on their critical line. For example, Figure~\ref{fig:rank 6}
depicts the $L$-function associated to the first known rank 6 elliptic
curve $E_6=[1,1,0,-2582,48720]$ of conductor $N=5\,187\,563\,742$ on the critical line,
and compares it to an elliptic curve of rank 0.

We are plotting the Hardy-$Z$ function, which is just $L_E(s)$ rotated so as to be
real on that line:
\begin{equation}
    \label{eq:Hardy-Z_E}
    Z_E(t) = w_E^{-1/2}
    \left(\frac{\sqrt{N}}{2\pi}\right)^{it} \exp(i \arg \Gamma(1+it)) L_E(1/2+it).
\end{equation}
One should be immediately struck by the large spikes that appear,
the first two occurring near $t=14$ and $t=21$. We compare this plot to
a less extreme $L$-function, say the elliptic curve of conductor 15, $E=[1,1,1,0,0]$,
also shown in the figure. $L$-functions tend to roll along, pulled by the rotation of
their Gamma factors and, here, the $(\sqrt{N}/(2\pi))^s$ factor that appears
in the functional equation.
These spikes are indeed very unusual, akin to finding Mount
Everest in Kansas. As such, they offer a clue into the phenomenon of rank in an
elliptic curve and should not be ignored.

The other obvious difference between the two plots, the higher frequency of the
$L$-function of $E_6$ compared to that of the conductor 15 curve, is easily
explained by its larger conductor, $5\,187\,563\,742$, which causes the
$(\sqrt{N}/(2\pi))^s$ factor to rotate faster and results
in an increase in the density of zeros of $L_{E_6}(s)$.

One quickly notices that the spikes occur when $t$ is near a zero $\gamma$
of the Riemann zeta function. For example, the first pair of zeros $1/2\pm i \gamma$
of the Riemann zeta function has $\gamma = 14.1347\ldots$, while the second has
$\gamma=21.022\ldots$. Because $L_E(1/2+it)$ gets large near these values, one can
surmise that the reciprocal of the Riemann zeta function must be
involved.

How can this arise?

\newpage

\newpage
\subsection{A bias in the $a(p)$'s}

Because the $L$-function is governed by its Dirichlet coefficients,
we should examine the $a(p)$'s. Table~\ref{table:ap e6} lists
the Dirichlet coefficients $a(p)$ of the rank 6 curve $E_6$, for $p\leq 173$.
\begin{table}[h!]
\centerline{
\begin{tabular}{|c|c||c|c||c|c||c|c|}
\hline
$p$ & $a(p)$ & $p$ & $a(p)$ & $p$ & $a(p)$ & $p$ & $a(p)$ \cr
\hline
2 & -1 & 31 & -8 & 73 & -11 & 127 & -20 \cr
3 & -1 & 37 & -11 & 79 & -10 & 131 & -20 \cr
5 & -4 & 41 & -10 & 83 & -8 & 137 & -9 \cr
7 & -4 & 43 & -11 & 89 & -6 & 139 & -12 \cr
11 & -6 & 47 & -10 & 97 & -10 & 149 & -6 \cr
13 & -6 & 53 & -13 & 101 & -7 & 151 & -10 \cr
17 & -7 & 59 & -3 & 103 & 4 & 157 & -9 \cr
19 & -8 & 61 & -10 & 107 & 4 & 163 & -14 \cr
23 & -8 & 67 & -12 & 109 & -13 & 167 & 12 \cr
29 & -6 & 71 & -15 & 113 & 0 & 173 & 4 \cr
\hline
\end{tabular}
}
\caption{The coefficients $a(p)$ of $E_6$, $p\leq 173$}.
\label{table:ap e6}
\end{table}

Notice that the initial $a(p)$'s are negative and, apart from the first four and $p=59$,
are less than or equal to $-6$ all the way through $p=101$. Indeed,
by Hasse's bound on $a(p)$, we can write
\begin{equation}
    \label{eq:a_p theta_p}
    a(p)/p^{1/2} = 2 \cos(\theta_p) = \alpha_p+\beta_p,
\end{equation}
where, for $p \nmid N$, $\alpha_p= e^{i\theta_p}$, and $\beta_p=e^{-i\theta_p}$, $\theta_p \in [0,\pi]$, and for $p \mid N$, $\alpha_p = a(p)$, $\beta_p = 0$.
Thus,
\begin{equation}
    \label{eq:L_p factors}
    {\mathcal{L}}_p(1/p^s) = (1-\alpha_p/p^s)^{-1}(1-\beta_p/p^s)^{-1},
\end{equation}
so that, taking the logarithmic derivative,
\begin{equation}
    \label{eq:log diff}
    \frac{L'(s)}{L(s)} = - \sum_{n=1}^\infty \frac{c(n)}{n^s},
\end{equation}
where
\begin{equation}
     \label{eq:c}
     c(n) = \begin{cases}
          \log(p) (\alpha_p^k+\beta_p^k), \quad &\text{if $n=p^k$,} \\
          0, \quad &\text{otherwise.}
     \end{cases}
\end{equation}

The Riemann and von-Mangoldt explicit formula for the coefficients $c(n)$ is thus
\begin{equation}
    \label{eq:explicit formula c(n)}
    \sum_{n\leq x} c(n) = -\sum_{\rho} \frac{x^\rho}{\rho} + o_E(x^{1/2}),
\end{equation}
where the sum is over the non-trivial zeros of $\rho$ of $L_E(s)$.
Now,
\begin{equation}
    \sum_{n\leq x} c(n) =
    \sum_{p\leq x} c(p)
    + \sum_{p\leq x^{1/2}} c(p^2) +O_{\varepsilon,E}(x^{1/3+\varepsilon}).
\end{equation}
Furthermore $c(p) = \log(p) (\alpha(p)+\beta(p)) = \log(p) a(p)/p^{1/2}$, and
$c(p^2) = \log(p) (\alpha(p)^2+\beta(p)^2) = \log(p) (a(p)^2-2p)/p$. One can
use Perron's formula on the symmetric square $L$-function of $L_E$ to show that
\begin{equation}
    \label{eq:sum c(p^2)}
    \sum_{p\leq x^{1/2}} c(p^2) \sim - x^{1/2}.
\end{equation}
Hence, if $L_E(s)$ has rank $r$, so that the term $\rho=1/2$ is counted $r$ times
in~\eqref{eq:explicit formula c(n)}, we have
\begin{equation}
    \label{eq:sum c(n) b}
    \sum_{p\leq x} \frac{\log(p)a(p)}{p^{1/2}}= -(2r-1) x^{1/2}
    -\sum_{\rho\neq 1/2} \frac{x^\rho}{\rho} + o_E(x^{1/2}).
\end{equation}
Letting
\begin{equation}
    \label{eq:average a_p}
    S_E(x) := \sum_{p\leq x} \log(p) a(p) = \sum_{p\leq x} \frac{\log(p) a(p)}{p^{1/2}} p^{1/2},
\end{equation}
summing by parts using~\eqref{eq:sum c(n) b}, and assuming the Riemann Hypothesis
for $L_E(s)$, the `bias' in $a(p)$ equals $-r+1/2$
in the sense of mean using the logarithmic density ~\cite{RS}~\cite{S}:
\begin{equation}
    \label{eq:log density mean}
    \frac{1}{\log{X}} \int_1^X \frac{S_E(x)}{x} \frac{dx}{x} = -r+1/2 + O_E(1/\log{X}).
\end{equation}
%Taking into account the non-trivial zeros of $L_E(s)$, we have, more precisely:
%\begin{eqnarray}
%    \label{eq:log density mean b}
%    &&\frac{1}{\log(X)/2} \int_{X^{1/2}}^X \frac{S_E(x)}{x} \frac{dx}{x} \notag \\
%    &&= -r+1/2 \notag \\
%    &&+\frac{2}{\log(X)/2}
%    \sum_{\gamma>0}
%    \frac{\cos(\gamma\log(X))-\cos(\gamma\log(X)/2)
%    - (\sin(\gamma\log(X))-\sin(\gamma\log(X)/2))/\gamma}{1+\gamma^2} \notag \\
%    &&\qquad +o_E(1/\log(X)),
%\end{eqnarray}
%where $1/2 \pm i\gamma$ denotes a typical conjugate pair of non-trivial zero of $L_E(s)$.

% \begin{equation}
%     \label{eq:sum c(n) c}
%     x^{-1/2} \sum_{p\leq x} \frac{\log(p) a(p)}{p^{1/2}}
% \end{equation}
% has mean $-2r+1$ on a logarithmic scale. See Sarnak's letter for an analysis of
% the bias in the distribution of the $a(p)$'s caused
% caused by the term $(2r-1)x^{1/2}$~\cite{S}, and the related paper of Fiorilli~\cite{F}.
% 
% However, by the prime number theorem,
% \begin{equation}
%     \label{eq:pnt}
%     \sum_{p\leq x} \frac{\log(p)}{p^{1/2}} \sim 2 x^{1/2}.
% \end{equation}
% 
% 
% that
% \begin{equation}
%     \label{eq:sum c(n) d}
%     x^{-1/2} \sum_{p\leq x} \frac{\log(p) (a(p)+r-1/2)}{p^{1/2}}
% \end{equation}
% has mean 0, on a logarithmic scale.

Thus, whenever an elliptic curve has rank $r$,
its $a(p)$'s are biased to the negative  by an amount equal to $-r$ in comparison
to the generic elliptic curve of rank 0, which has a bias in the $a(p)$'s just of size
$1/2$.

We list the values of the numerically computed bias, i.e. the lhs
of~\eqref{eq:log density mean}, with $X=10^8$, for the curves $E_1$, $E_2$,
$E_3$, $E_4$, $E_5$, $E_6$, $E_7$, $E_{11}$, $E_{24}$: $-0.63$, $-1.58$, $-2.49$,
$-3.36$, $-4.22$, $-5.09$, $-6.04$, $-9.19$, $-16.58$. These values are rounded, and are
reasonably close to $-r+1/2$, especially, here, for $r\leq 7$.
Taking into account the non-trivial
zeros of $L_E(s)$, the implied constant
in the $O_E$ term in~\eqref{eq:log density mean} is typically of size
$(\log{N})\log((\log{N})/r)$, though can be larger, depending on the
location of the zeros of $L_E(s)$ near the real axis. This explains why the
agreement is not as good, here, for $r=11,24$ where the conductors are quite
large.

\begin{figure}[H]
\centerline{
    \psfig{figure=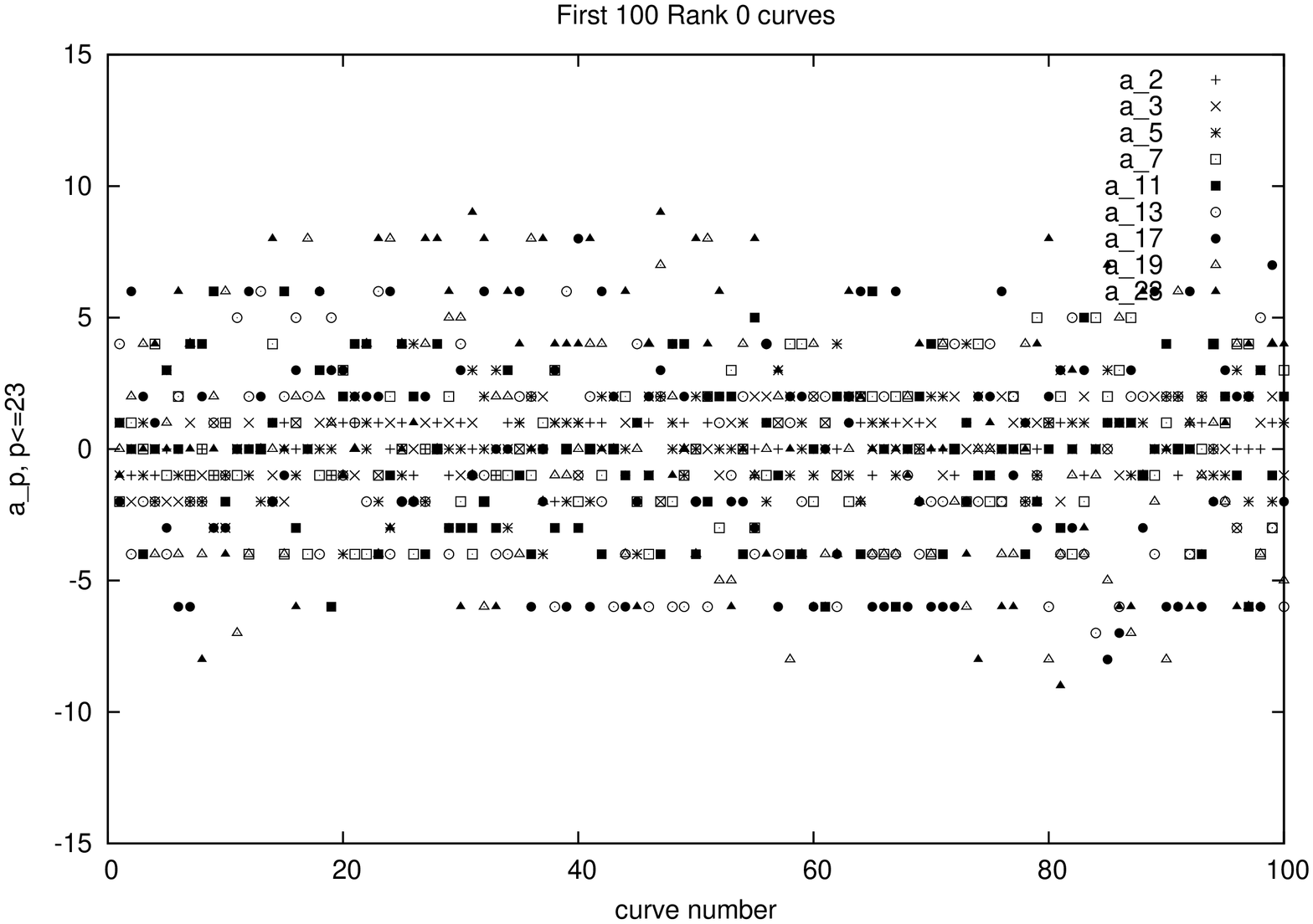,width=2.1in,height=3.5in,angle=-90}
    \psfig{figure=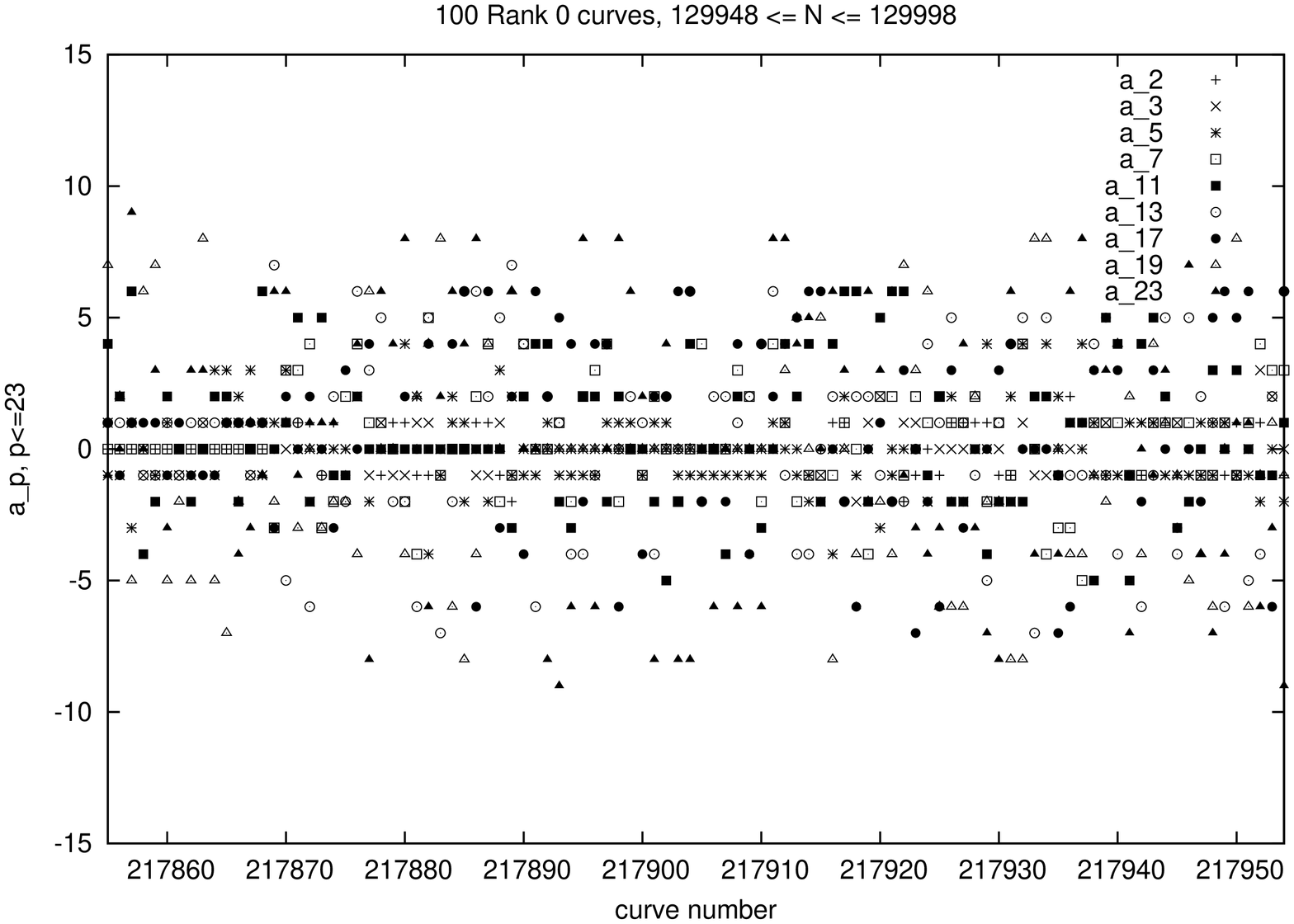,width=2.1in,height=3.5in,angle=-90}
}
\centerline{
    \psfig{figure=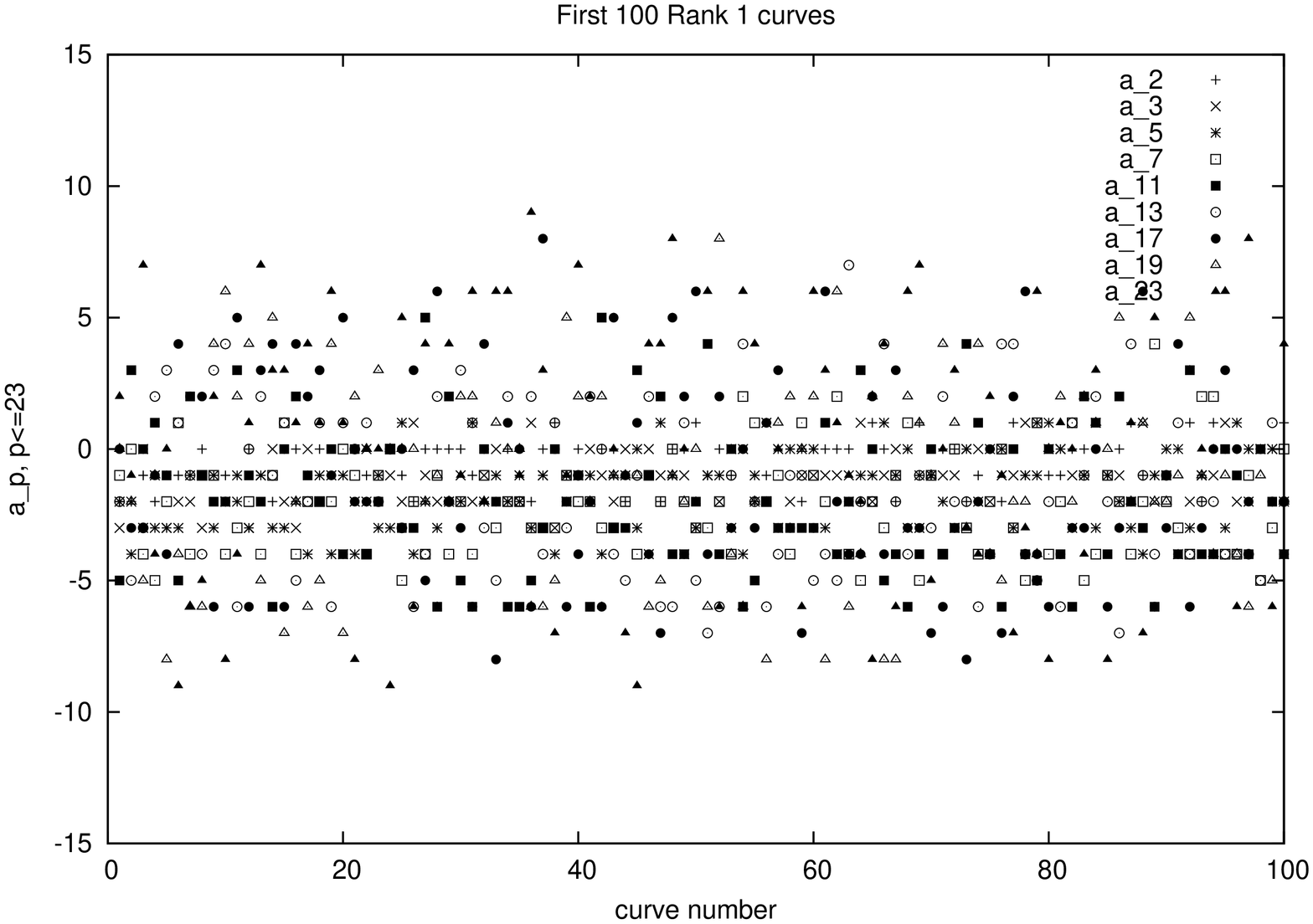,width=2.1in,height=3.5in,angle=-90}
    \psfig{figure=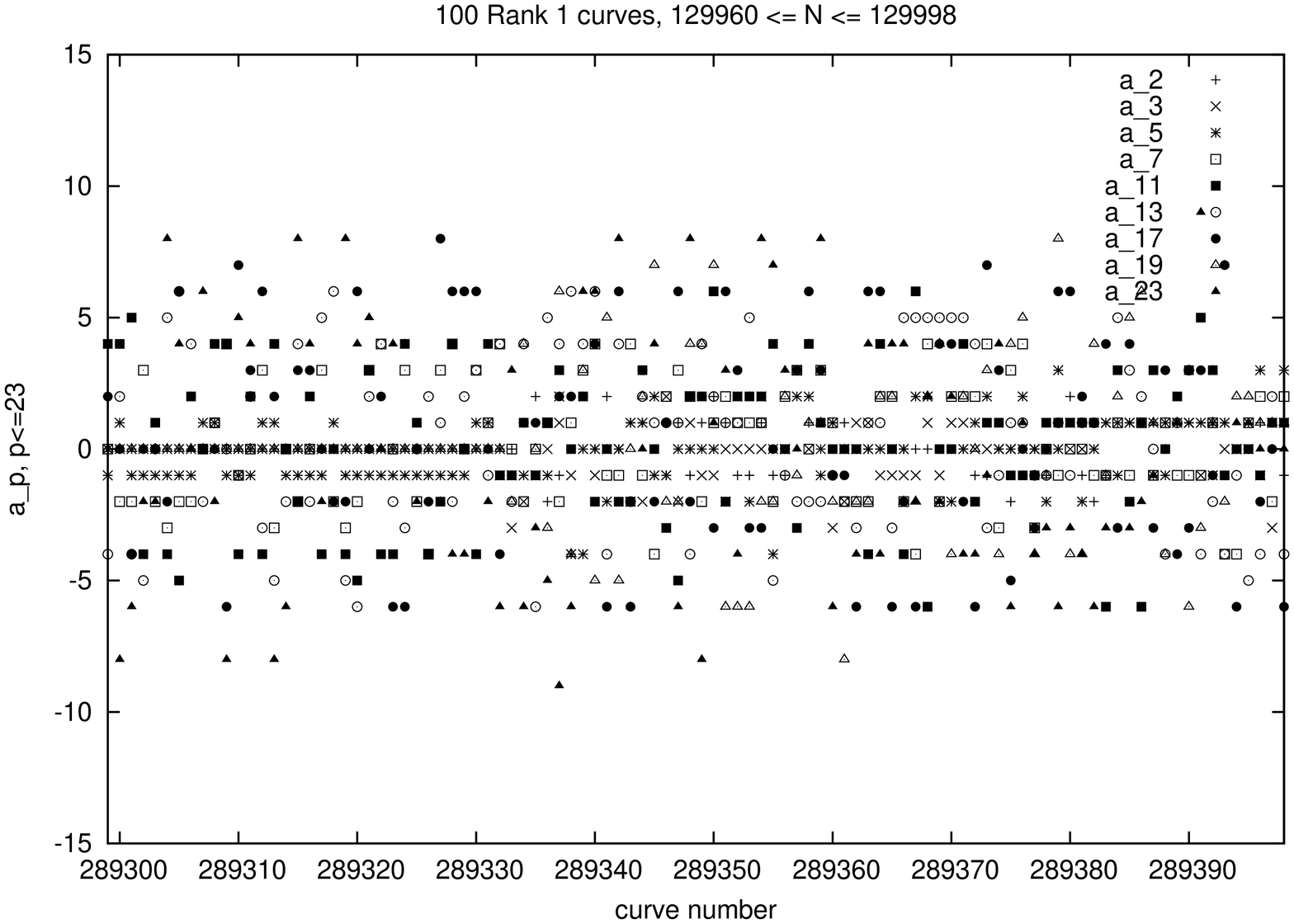,width=2.1in,height=3.5in,angle=-90}
}
\centerline{
    \psfig{figure=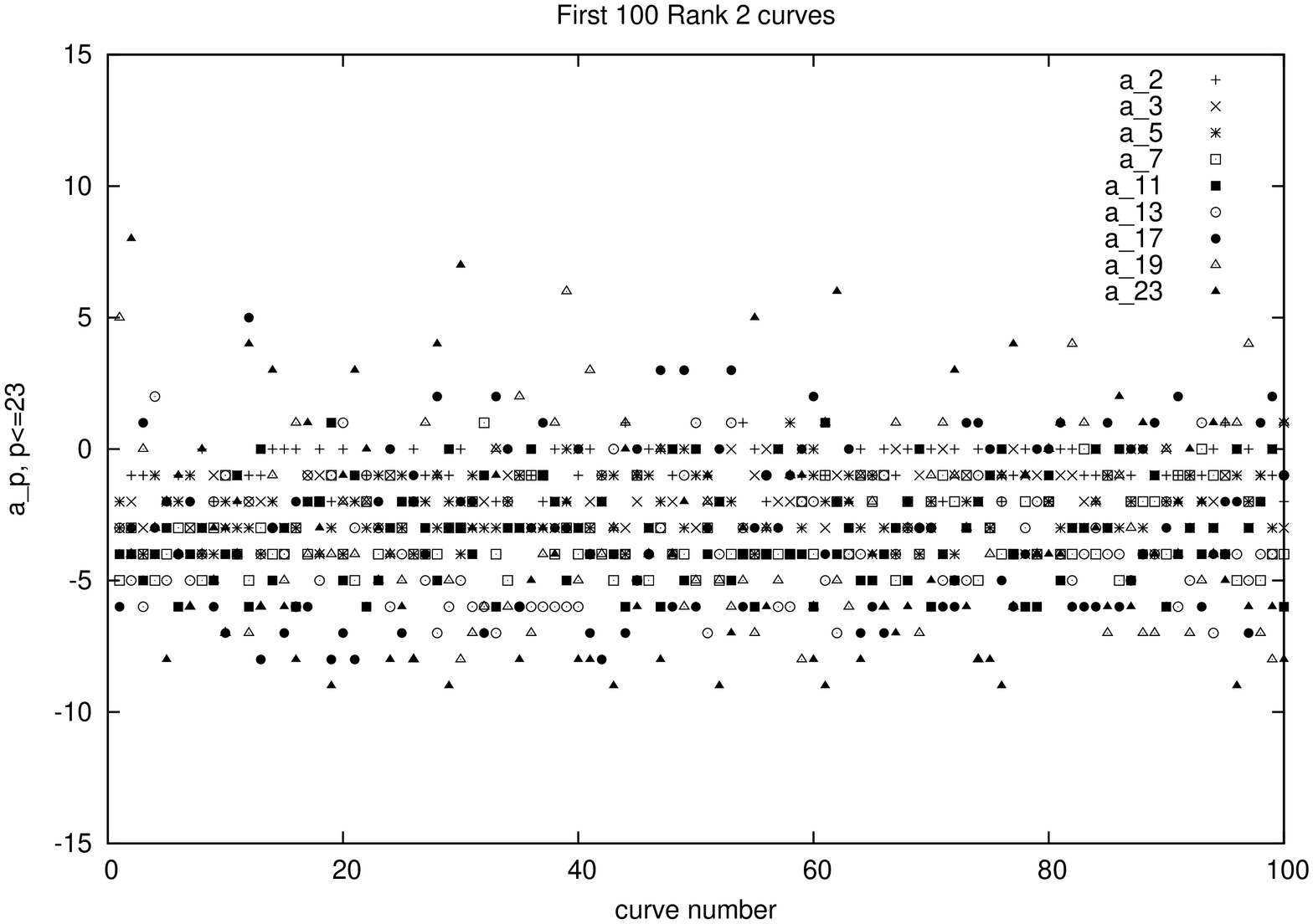,width=2.1in,height=3.5in,angle=-90}
    \psfig{figure=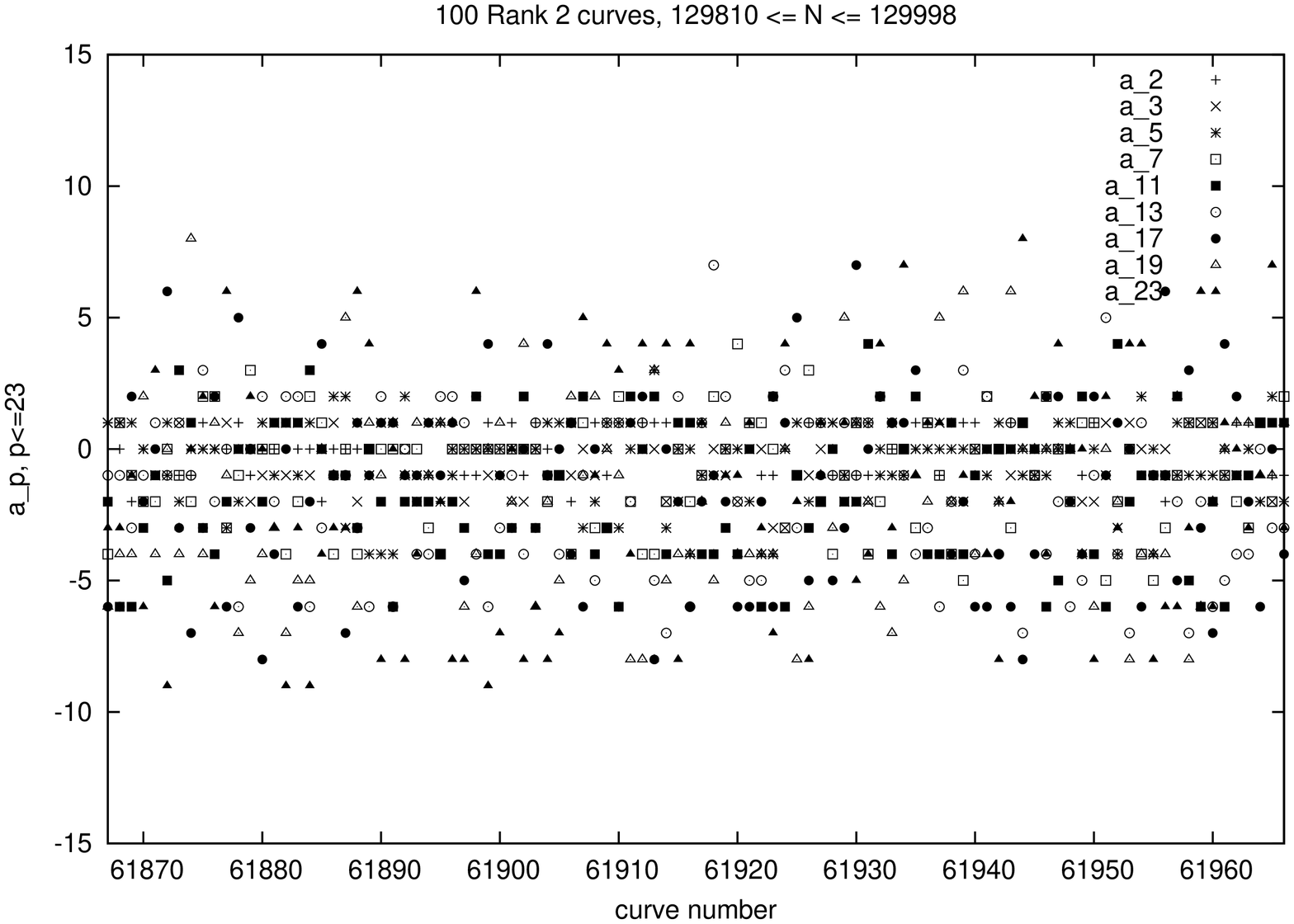,width=2.1in,height=3.5in,angle=-90}
}
\centerline{
    \psfig{figure=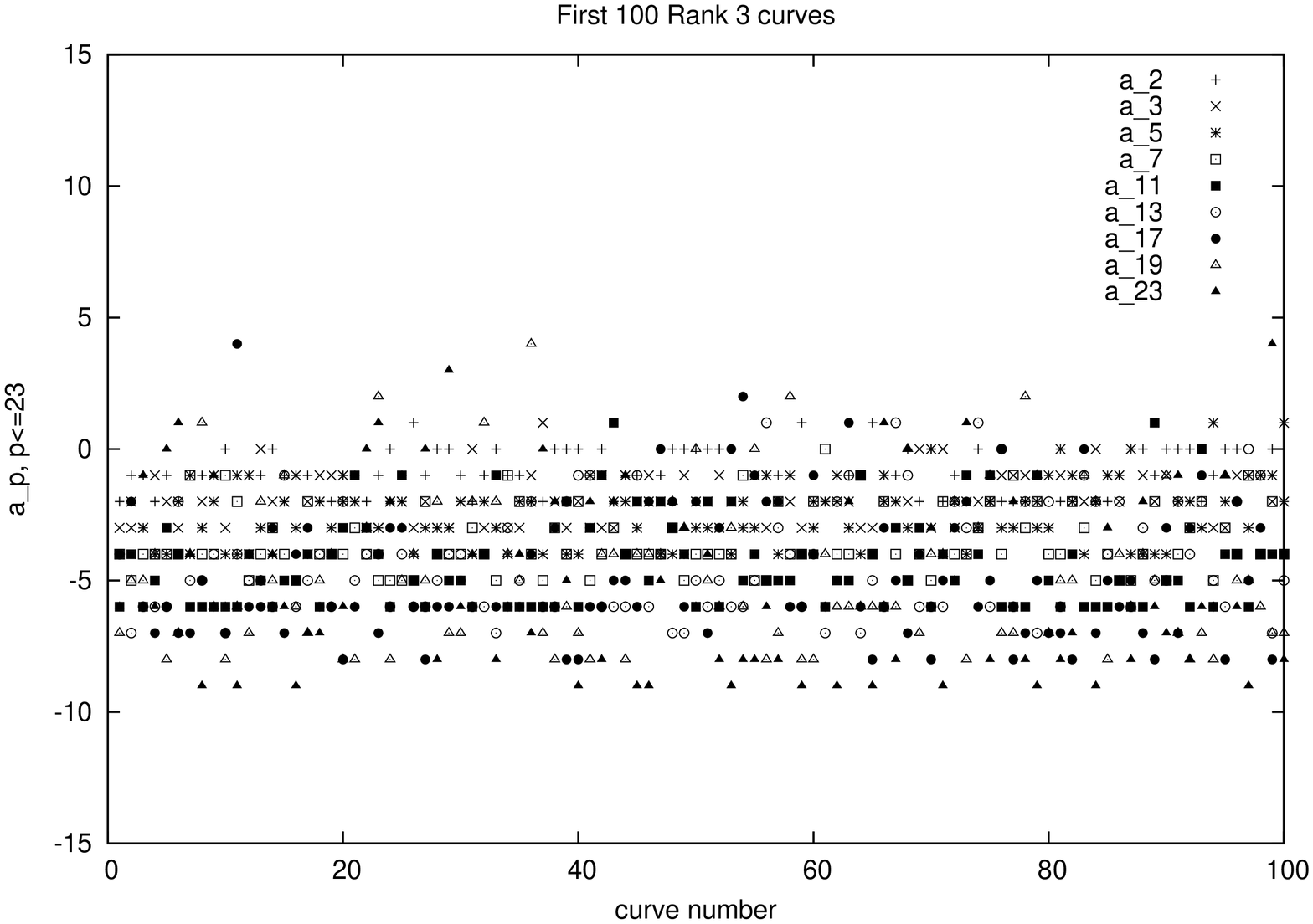,width=2.1in,height=3.5in,angle=-90}
    \psfig{figure=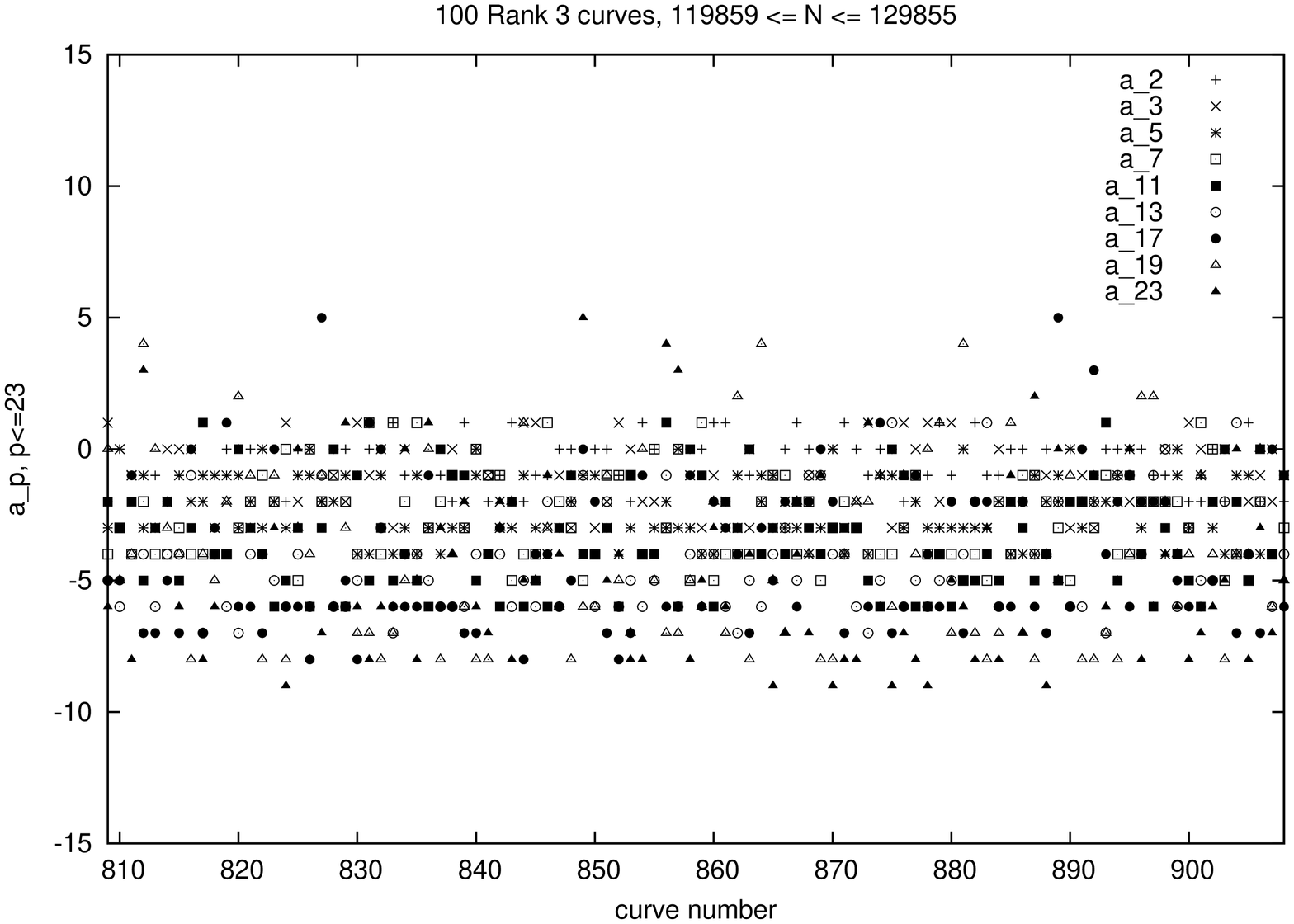,width=2.1in,height=3.5in,angle=-90}
}
\caption{Plots of $a(p)$, $p\leq 23$ for the first 100 elliptic curves of ranks
0, 1, 2,and 3, and also 100 curves, for each of these ranks, of conductor around
130,000. As the rank increases, the $a_p$'s tend to move downwards.}
\label{fig:a_p rank 0-3}
\end{figure}

We depict, in Figure~\ref{fig:a_p rank 0-3}, the values of $a(p)$, $p\leq 23$,
for the first 100 elliptic curves of ranks 0, 1, 2, and 3, as well as for 100
elliptic curves of these ranks but conductors of size approximate 130,000. The
curves were obtained from Cremona's database~\cite{C}. One can literally see
that, as a whole, the $a(p)$'s have a preference, on average to move downwards
as the rank increases.

We also graph the values of $p$ versus $a(p)$ for the rank 11 elliptic curve of smallest known
conductor, $E_{11}=[0, 0 ,1 ,-16359067 ,26274178986]$, found by Elkies and Watkins,
and for an elliptic curve of rank at least 24, $E_{24}$, of Martin and McMillen~\cite{MM}.
A large bias is very evident in these pictures, with most of the initial $a(p)$'s being negative.

\newpage
\newpage

\begin{figure}[H]
\centerline{
    \psfig{figure=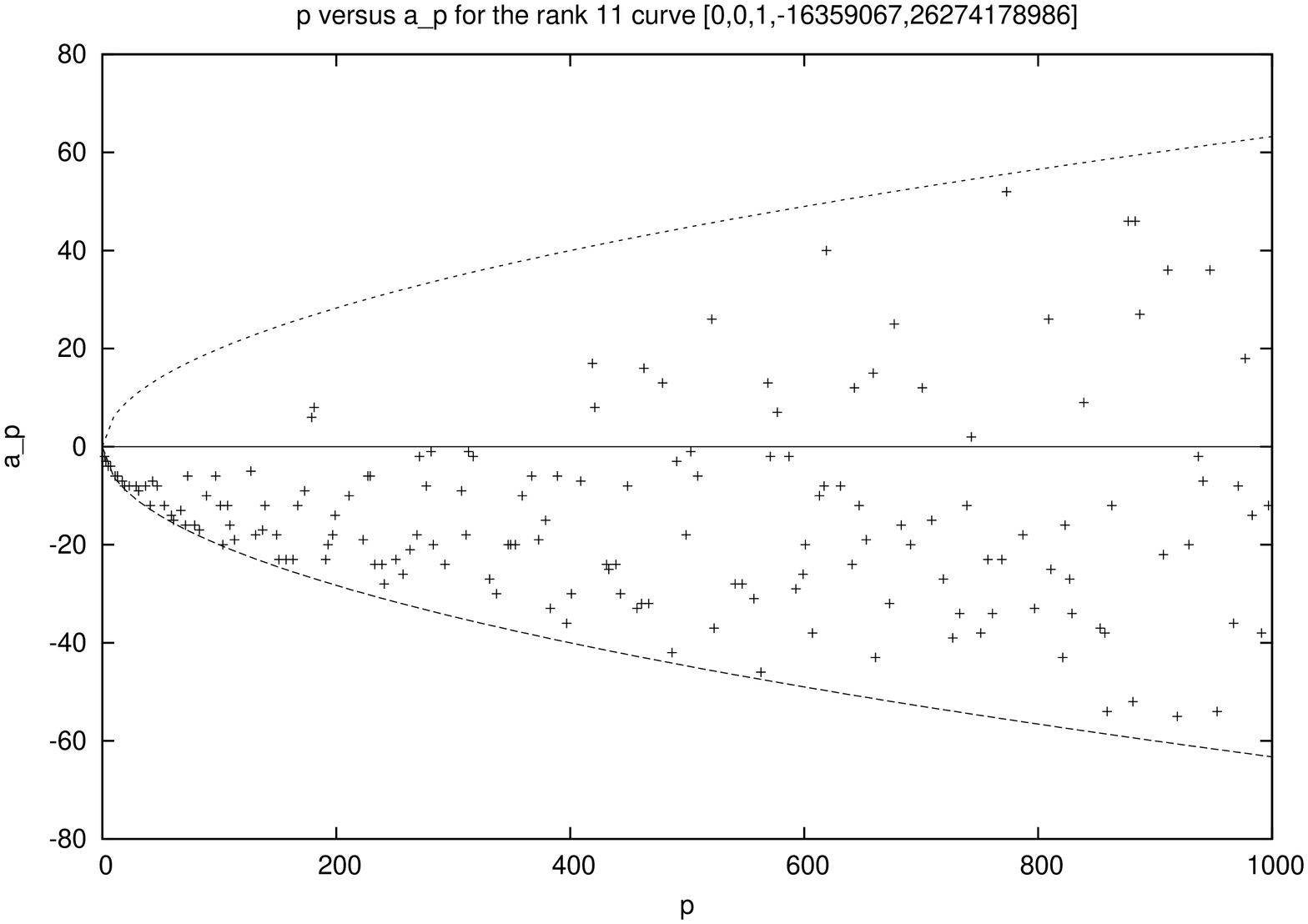,width=2.5in,angle=-90}
    \psfig{figure=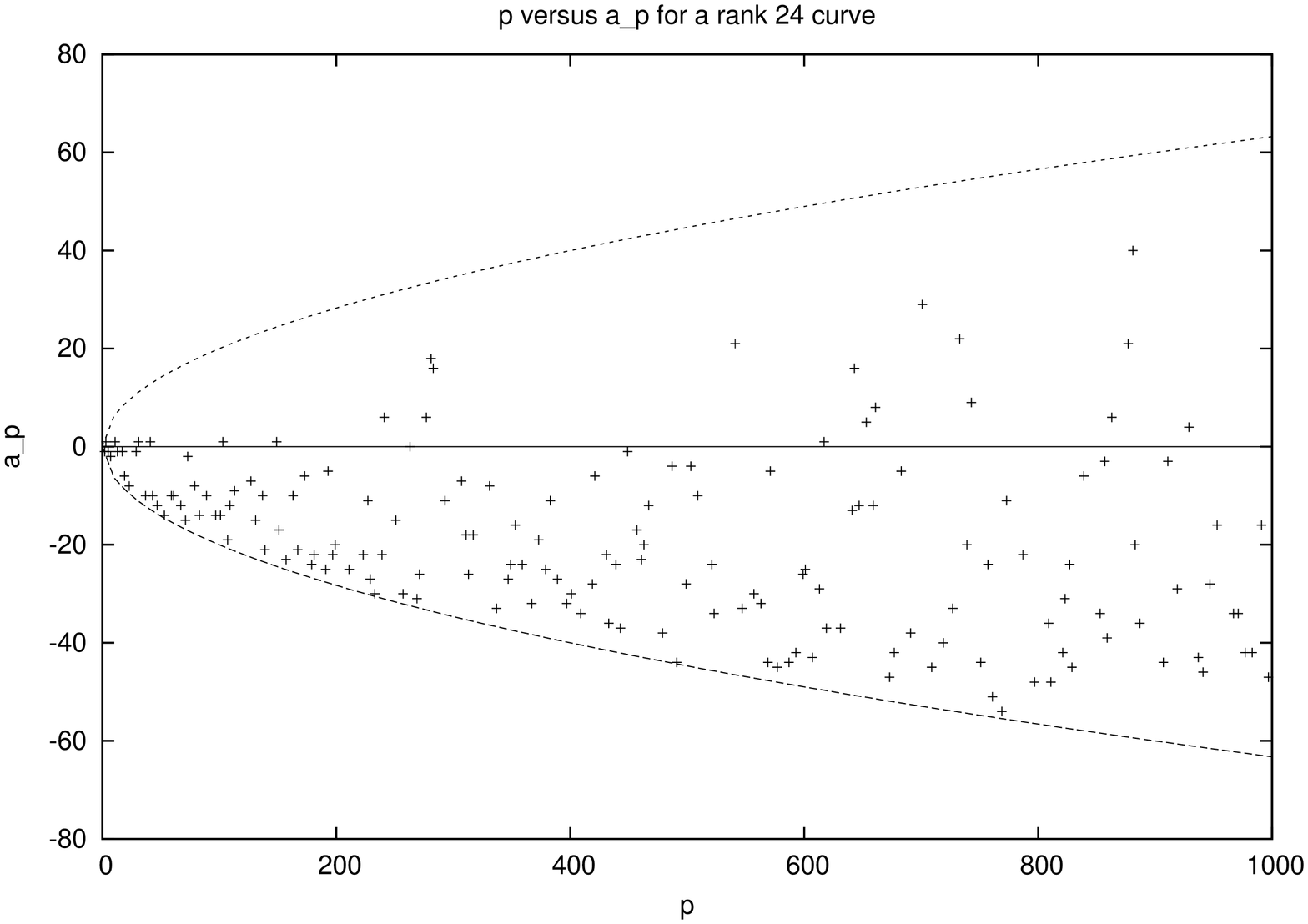,width=2.5in,angle=-90}
}
\caption{Plot of $a(p)$, $p\leq 1000$ for $E_{11}$, the rank 11 curve of smallest known conductor,
and for $E_{24}$ of rank 24. The dashed curve is Hasse's bound $\pm 2 p^{1/2}$.}
\label{fig:a_p ranks 11,24}
\end{figure}

For the curve of rank 11, where the conductor is comparatively small relative to the
rank, we notice that many of the initial $a(p)$'s are very close to the negative
side of Hasse's bound, $-2p^{1/2}$. It would be worthwhile to analyse
this feature via Weil's explicit formula, but we do not consider that
in this paper.

\subsection{Accounting for the bias in the $a(p)$'s}

The typical factor appearing in~\eqref{eq:euler_prdct} can be expanded:
\begin{equation}
    \label{eq:local factor expanded}
    \left(1-a(p) p^{-(s+1/2)}+p^{-2s}\right)^{-1}
    = 1 + \frac{a(p)}{p^{1/2}}\frac{1}{p^s}+\frac{(a(p)^2-p)}{p} \frac{1}{p^{2s}} + \ldots.
\end{equation}
Consider, now, the elliptic curve of rank 6, $E_6$, which has bias in $a(p)$ of size $-6+1/2$.
Because an elliptic curve of rank 0 has a bias that is equal to $1/2$, 
we experimented by pulling out, from each local factor in the Euler product for $L_{E_6}(s)$,
a factor of
\begin{equation}
    \left(1- p^{-s-1/2}\right)^6 = 1 - 6 p^{-s-1/2} + \ldots,
\end{equation}
in order to account for the excess bias in $a(p)$ of size $-6$.
We thus write
\begin{equation}
    L_{E_6}(s) = \frac{1}{\zeta(s+1/2)^6} f_{E_6}(s),
\end{equation}
and then test whether $f_{E_6}(s):=L_{E_6}(s)/\zeta(s+1/2)^6$ is relatively
tempered compared to $L_{E_6}(s)$.

The factor of $1/\zeta(s+1/2)^6$ not only captures terms in the Euler product, but
also the extreme behaviour of $L_{E_6}(s)$, both its large spikes and it's 6th
order vanishing at $s=1/2$.
The latter arises from the first order zero of $1/\zeta(s+1/2)$ at $s=1/2$.

\begin{figure}[H]
    \centerline{
            \psfig{figure=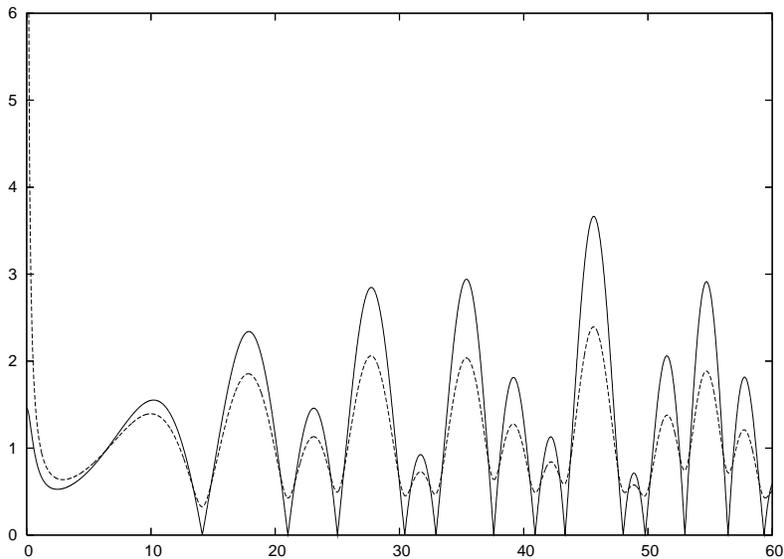,width=3in,angle=-90}
    }
    \caption
    {
        $|\zeta(1+it)|$ (dotted) compared to $|\zeta(1/2+it)|$ (solid).
    }
    \label{fig:zeta one line}
\end{figure}

The spikes of $L_{E_6}(1/2+it)$ near the zeros of $\zeta$ are also accounted
for by the factor $1/\zeta(s+1/2)^6$. Even though we are evaluating the zeta
function on the one line $s=1+it$, the minima, in $t$, of $|\zeta(1+it)|$ occur
near the minima $|\zeta(1/2+it)|$, at least initially, and hence
$1/\zeta(1+it)$ spikes when $t$ is near a zero of zeta.
This can be explained via the Hadamard product for $\zeta$. A related phenomenon
is discussed later in the paper, around equation~\eqref{eq:zeta log diff}.
Figure~\ref{fig:zeta one line}, taken from~\cite{R3}, compares $|\zeta(s)|$ on the
$1/2$ and $1$ lines,
and illustrating that the minima of both roughly coincide initially. See also
the third plot of Figure~\ref{fig:3 plots one line} which depicts $1/|\zeta(1+it)|$.

The left plot in
Figure~\ref{fig:spikes} confirms that the large spikes of $L_E(1/2+it)$ occur
near the maxima of $1/|\zeta(1+it)|^6$. We notice that the peaks are too high.
However, it seems that we can fix this by adjusting for the initial primes for
which Hasse's bound prevents us from achieving $a(p)= -r$. These small primes,
while few in number, have a dramatic effect on the graph of $L_{E_6}(1/2+it)$.
Therefore, we compare $L_{E_6}(1/2+it)$, in the right plot of
Figure~\ref{fig:spikes}, to the function
\begin{equation}
    \label{eq:with local correction}
    \frac{\text{local}_E(1/2+it)}{\zeta(1+it)^6}
\end{equation}
with
\begin{equation}
    \label{eq:local E6}
    \text{local}_{E_6}(s)= f_{E_6}(2,s)f_{E_6}(3,s)f_{E_6}(5,s)f_{E_6}(7,s)
\end{equation}
where $f_{E_6}(p,s)$ corrects the local factor of $1/\zeta(s+1/2)^6$, at the early primes $p<9$
for which Hasse's bound prevents $a(p)$ from achieving its bias of $-6$,
to match those of $L_{E_6}(s)$ in (\ref{eq:euler_prdct}):
In this example,
\begin{eqnarray}
    \notag
    f_{E_6}(2,s) &&= (1+2^{-s-1/2})^{-1} (1-2^{-s-1/2})^{-6}\notag \\
    f_{E_6}(3,s) &&= (1+3^{-s-1/2})^{-1} (1-3^{-s-1/2})^{-6} \notag \\
    f_{E_6}(5,s) &&= (1+ 4\cdot 5^{-s-1/2}+5^{-2s})^{-1} (1-5^{-s-1/2})^{-6}\notag \\
    f_{E_6}(7,s) &&= (1+ 4\cdot 7^{-s-1/2}+7^{-2s})^{-1} (1-7^{-s-1/2})^{-6}.
\end{eqnarray}

\afterpage{
\begin{figure}[H]
    \centerline{
            \psfig{figure=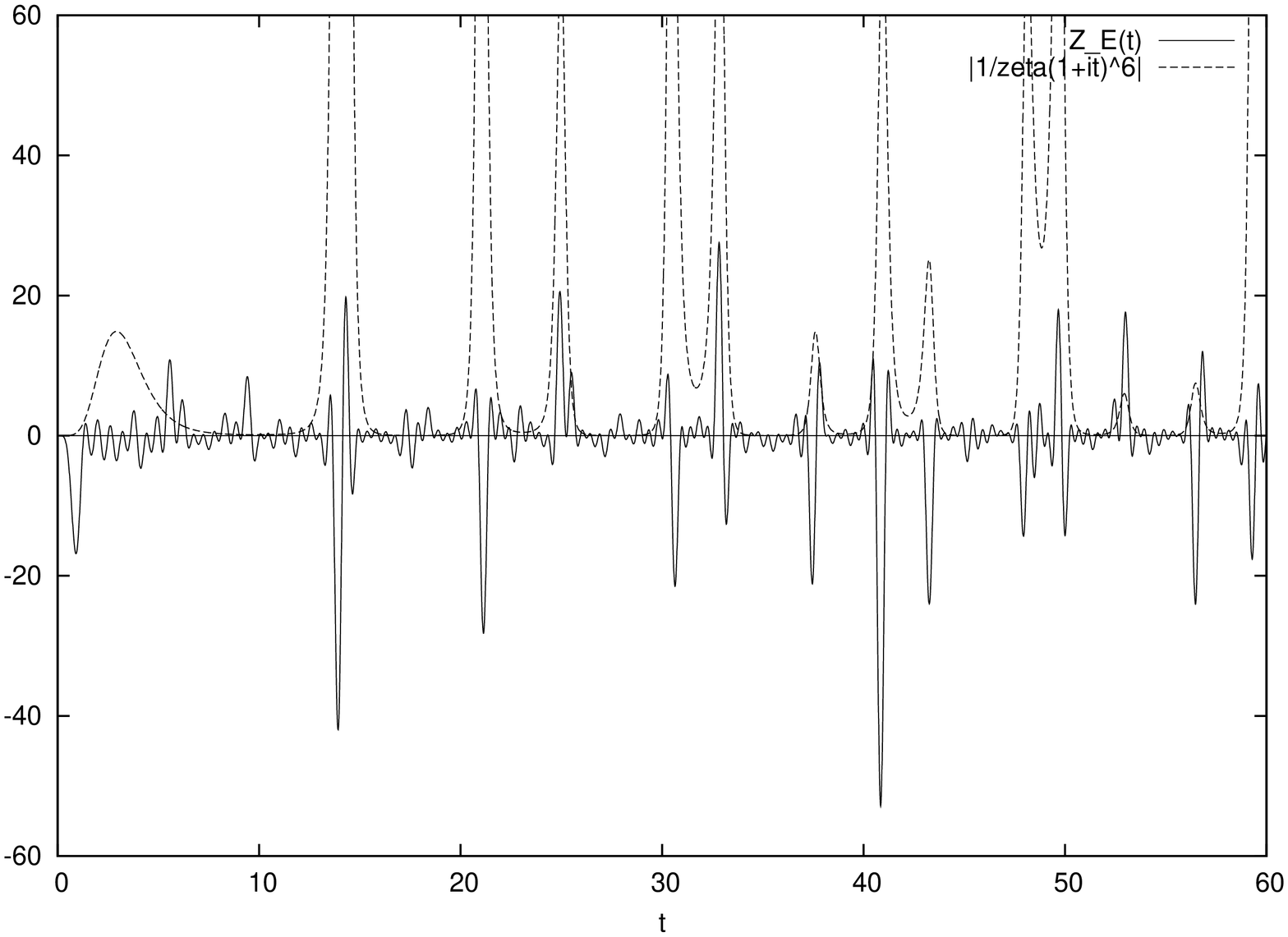,width=2.5in,angle=-90}
            \psfig{figure=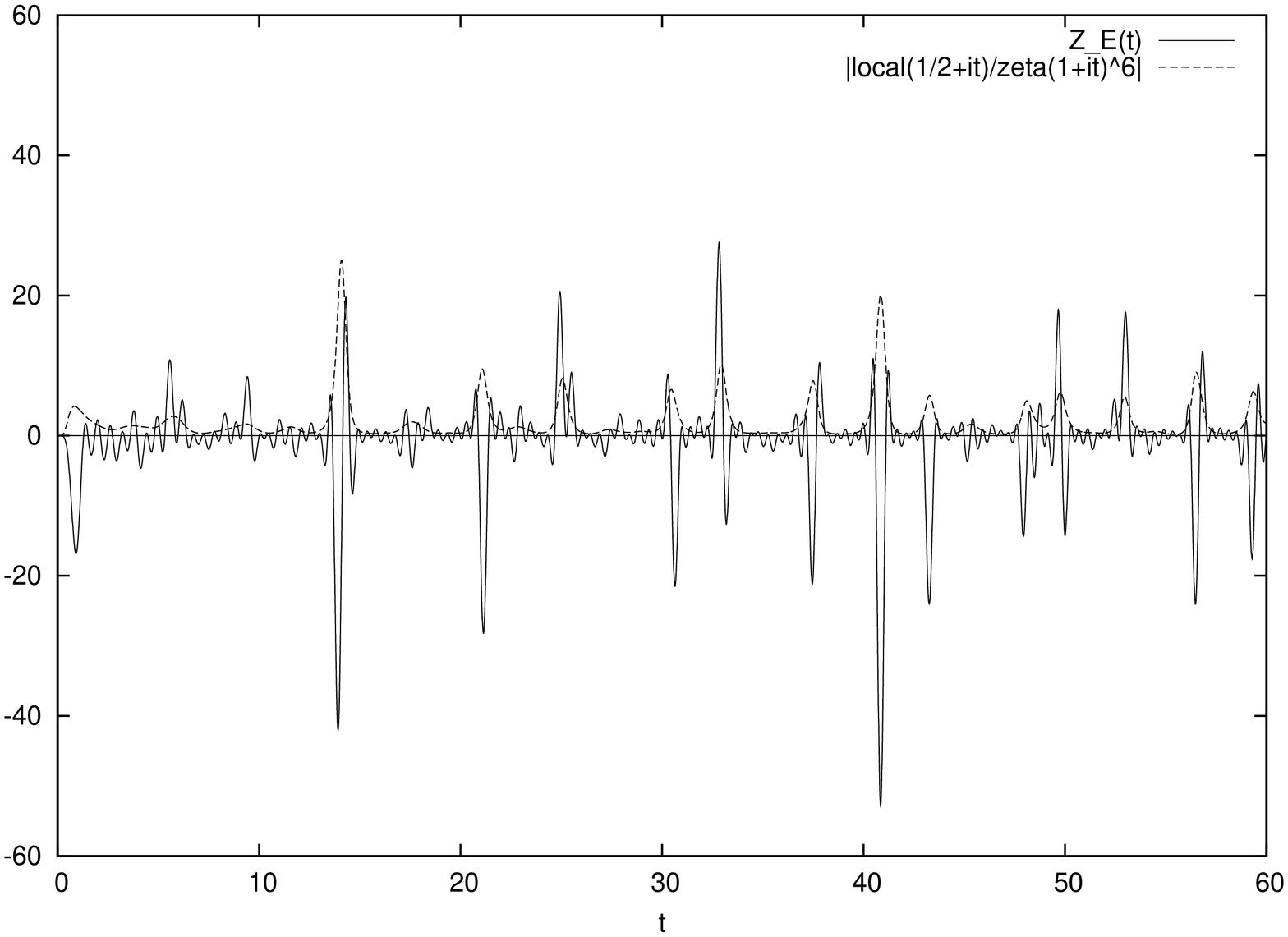,width=2.5in,angle=-90}
    }
    \caption
    {
         The Hardy function $Z_{E_6}(t)$ (solid) compared, in the left plot,
         to $1/|\zeta(1+it)|^6$ (dotted), and, in the right plot, to
         $|\text{local}_{E_6}(1/2+it)/\zeta(1+it)^6|$.
    }
    \label{fig:spikes}
\end{figure}
The following graph depicts the ratio $Z_{E_6}(t)/ (|\text{local}_{E_6}(1/2+it)/\zeta(1+it)^6|)$.
\begin{figure}[H]
    \centerline{
            \psfig{figure=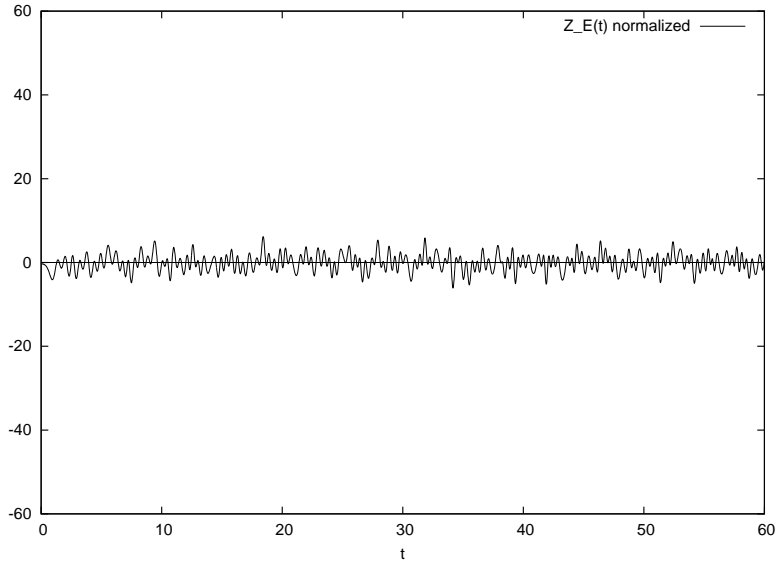,width=3in,angle=-90}
    }
    \caption
    {
        $Z_{E_6}(t) |\zeta(1+it)^6/\text{local}_{E_6}(1/2+it)|$
    }
    %\label{fig:conj vs reality}
\end{figure}
}

Other examples have the same features, with the peaks being more prominent when the rank is larger.
Below we display, for $0 \leq t \leq 100$, the graphs of $L_{E_r}(1+it)$, for
$r=1,2,3,4,5,7$.
We also define local correction functions for each of these seven elliptic curves
for the local factors for which $2p^{1/2}< r$.
For $E_1$ and $E_2$, $\text{local}_E(s)=1$.
For $r=3, 4, 5, 7$ it involves correcting for the primes $p\leq 2,3,5,11$ respectively.
These plots were generated using the author's $L$-function computer package, {\tt lcalc},
which uses a smooth approximate functional equations to compute $L$-functions. It also
relies on {\tt PARI}'s elliptic curve routines for computing the $a(p)$'s and conductor
associated to an elliptic curve~\cite{R2}~\cite{P}.

\afterpage{
\begin{figure}[H]
\centerline{
    \psfig{figure=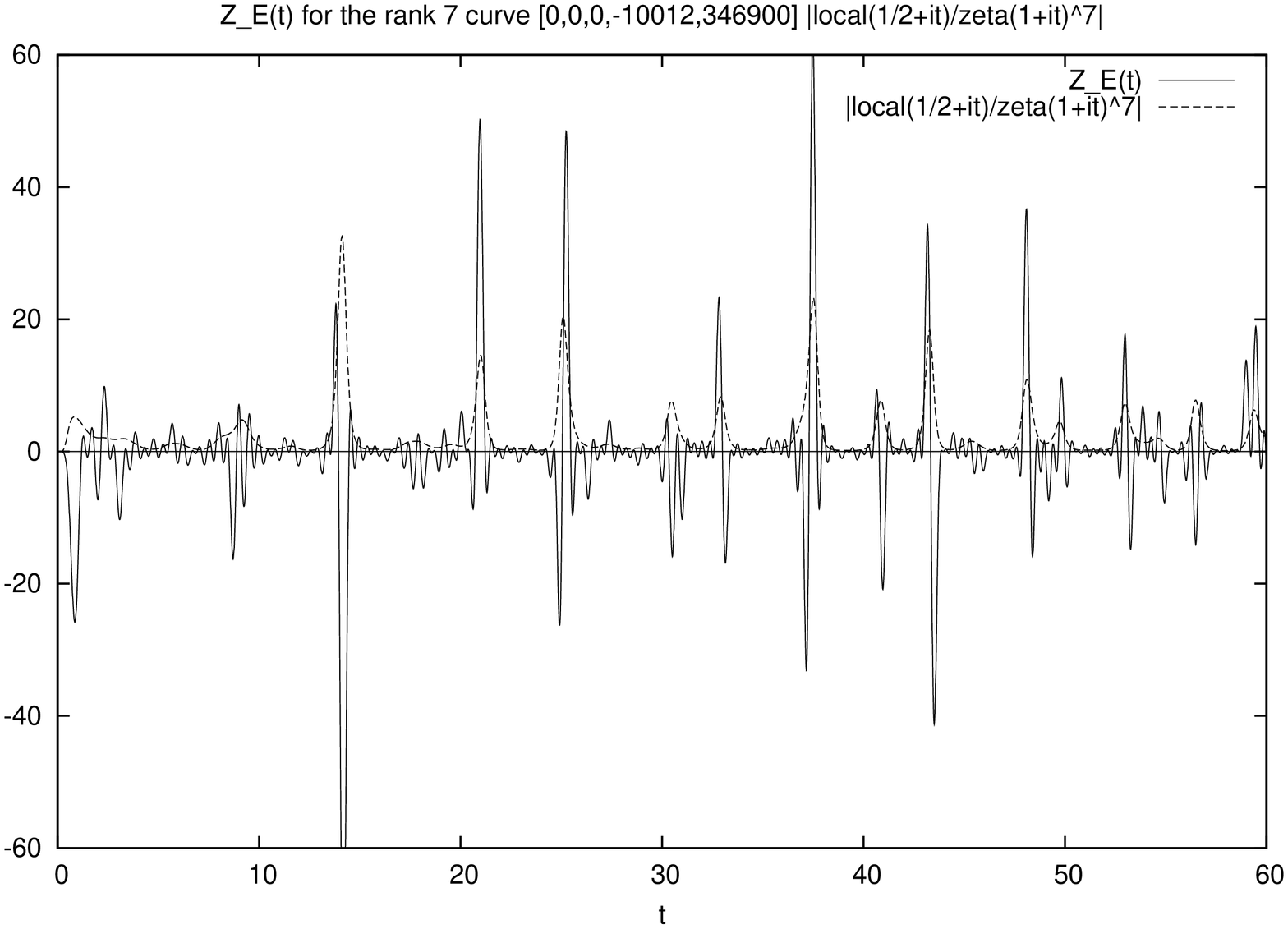,width=2.5in,angle=-90}
    \psfig{figure=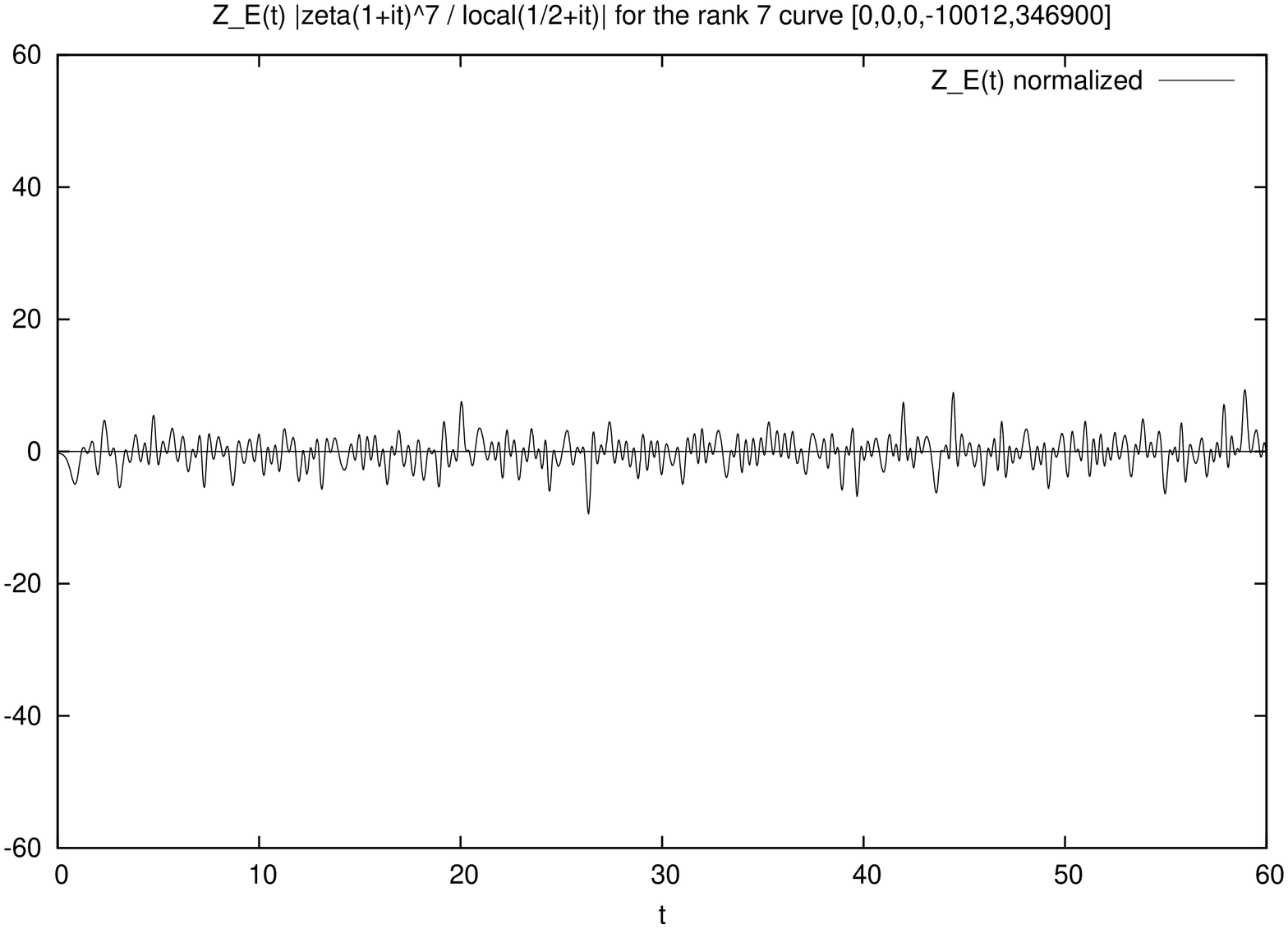,width=2.5in,angle=-90}
}
%\centerline{
%    \psfig{figure=e6_withzeta.ps,width=3.5in,angle=-90}
%    \psfig{figure=e6_normalized.ps,width=3.5in,angle=-90}
%}
\centerline{
    \psfig{figure=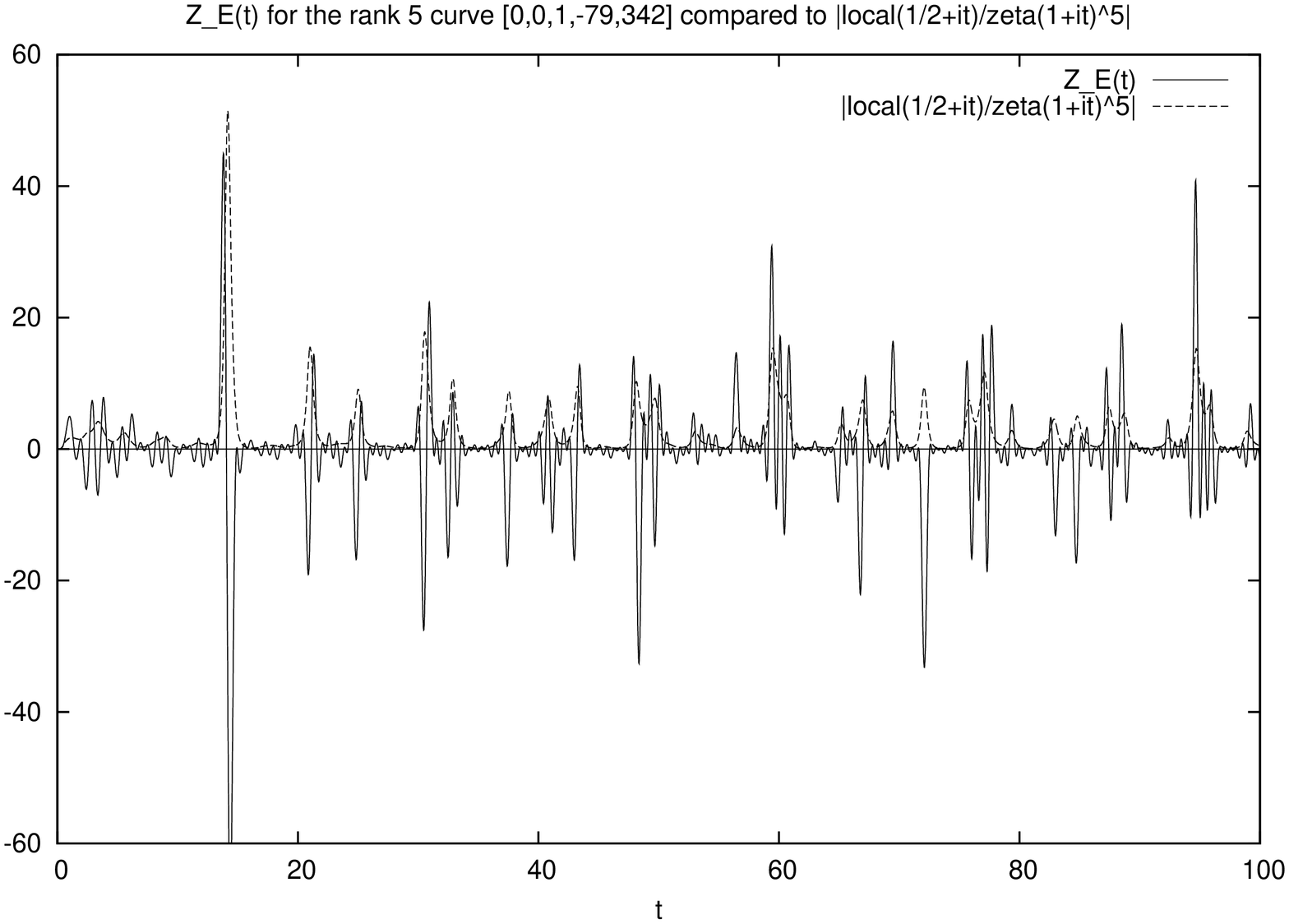,width=2.5in,angle=-90}
    \psfig{figure=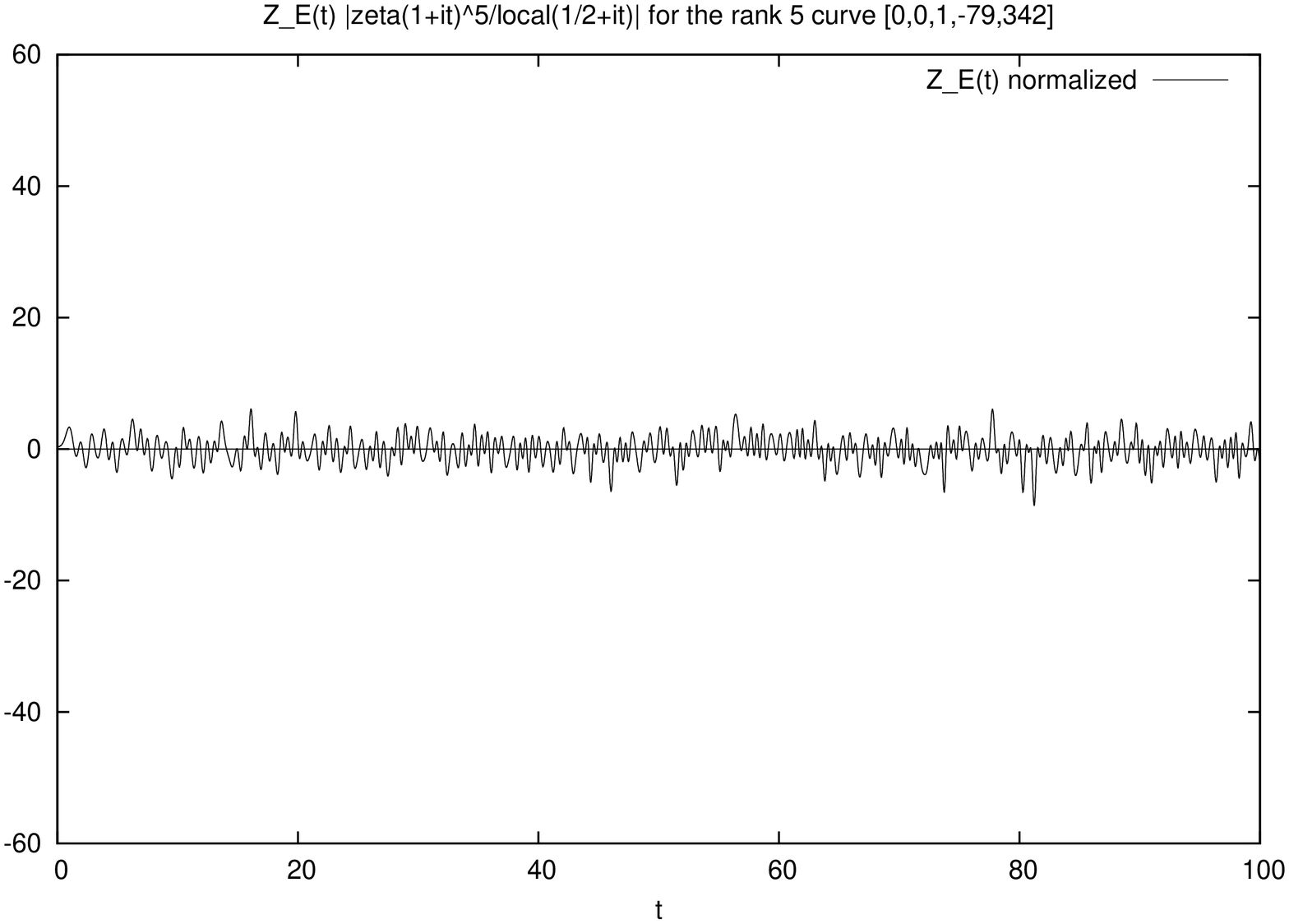,width=2.5in,angle=-90}
}
\centerline{
    \psfig{figure=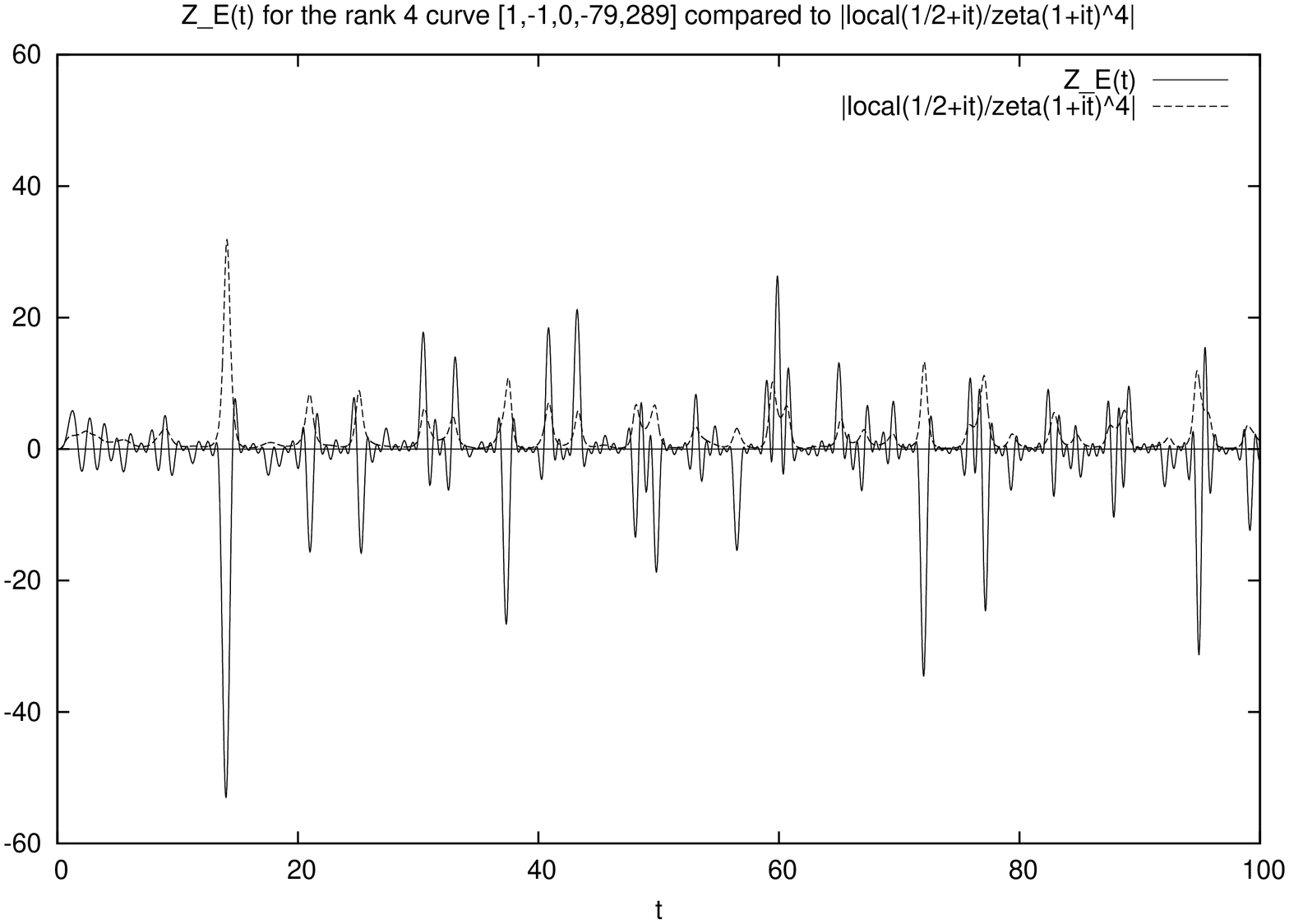,width=2.5in,angle=-90}
    \psfig{figure=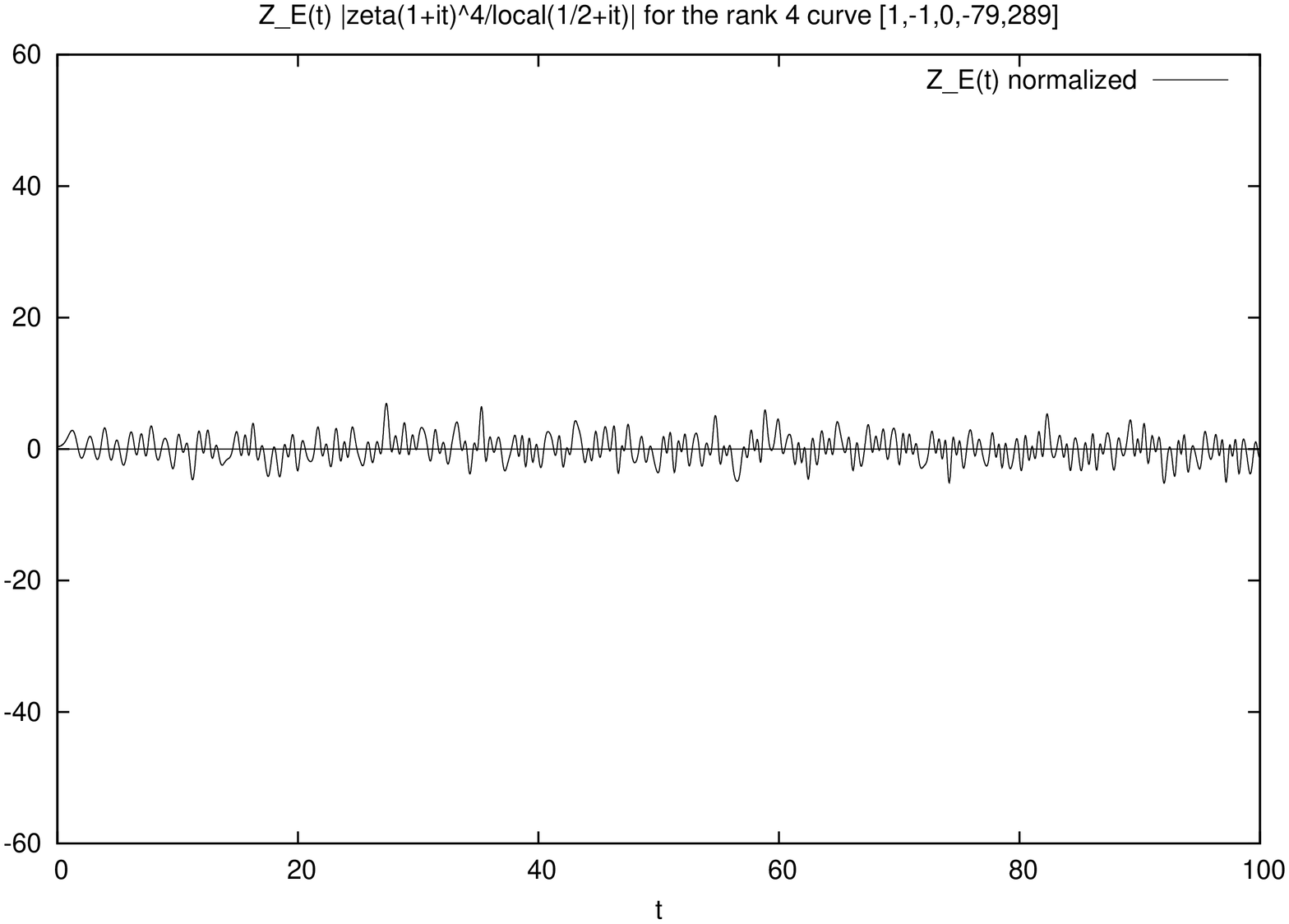,width=2.5in,angle=-90}
}
\caption
    {Graphs in the left column show $Z_{E_r}(t)$ (solid), for $r=7,5,4$, compared to
    $|\text{local}_{E_r}(1/2+it)/\zeta(1+it)^r|$ (dotted). Graphs in the right column
    show $Z_{E_r}(t) |\zeta(1+it)^r/\text{local}_{E_r}(1/2+it)|$.}
\label{fig:graphs 7-5}
\end{figure}
}

\afterpage{
\begin{figure}[H]
\centerline{
    \psfig{figure=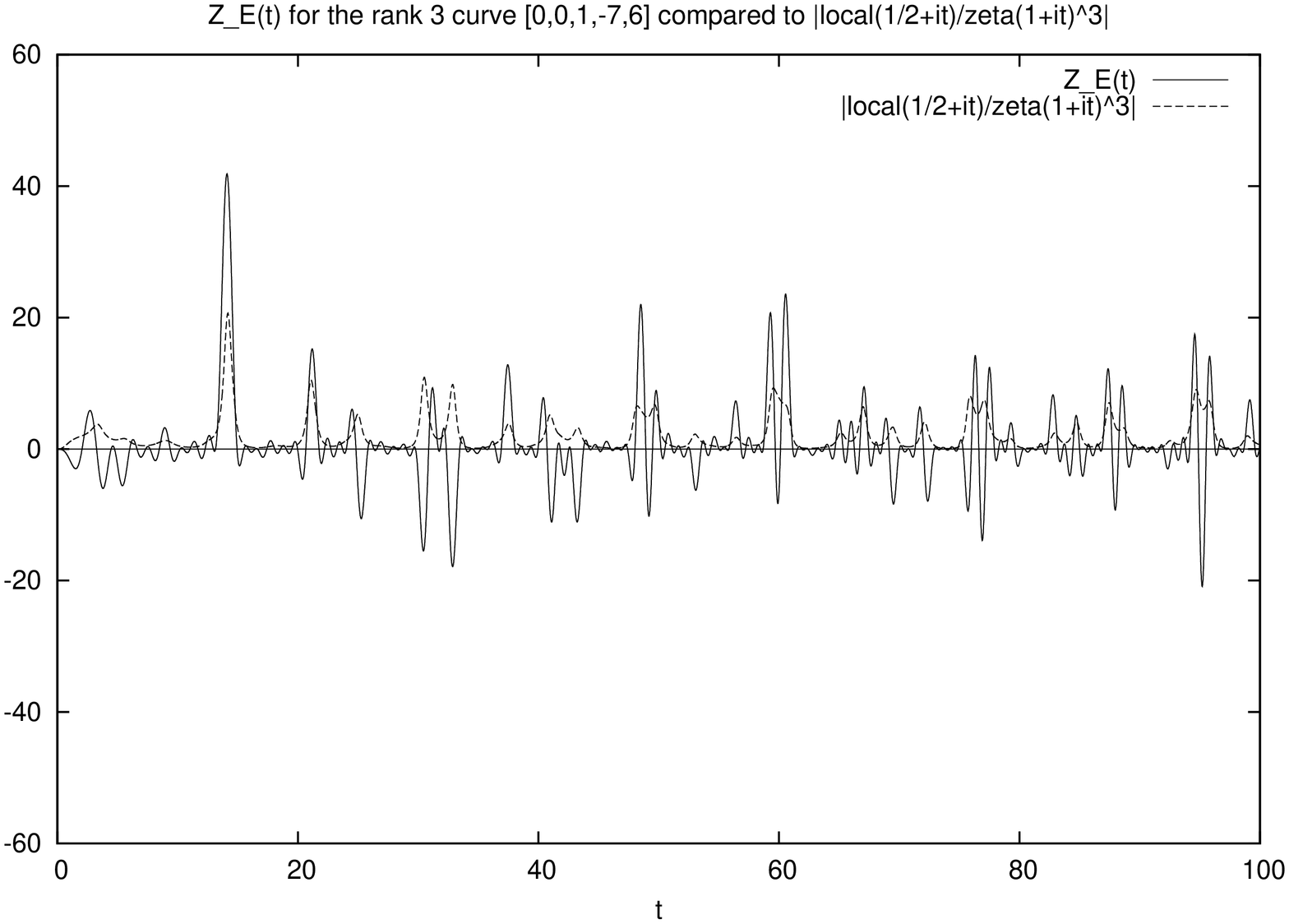,width=2.5in,angle=-90}
    \psfig{figure=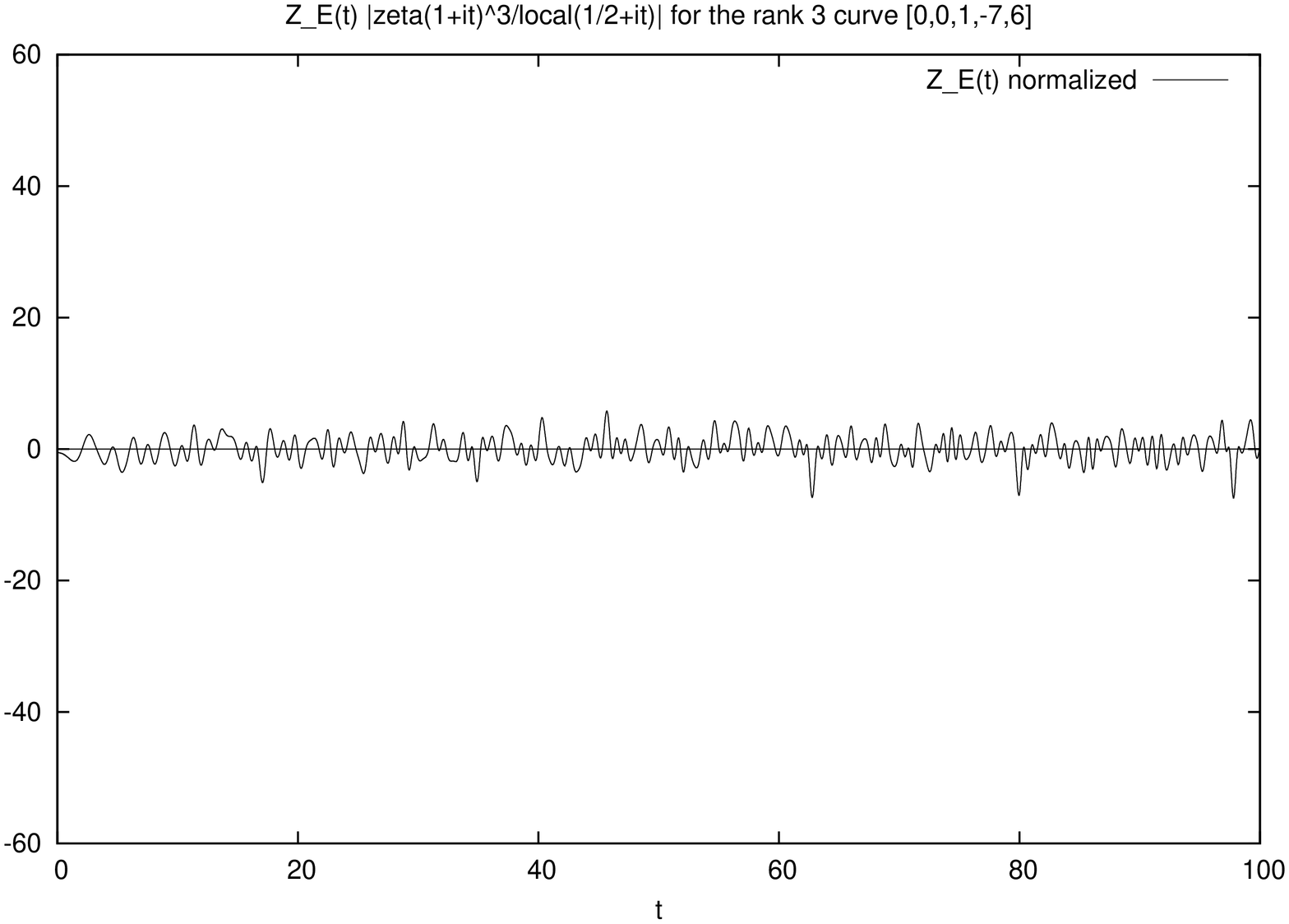,width=2.5in,angle=-90}
}
\centerline{
    \psfig{figure=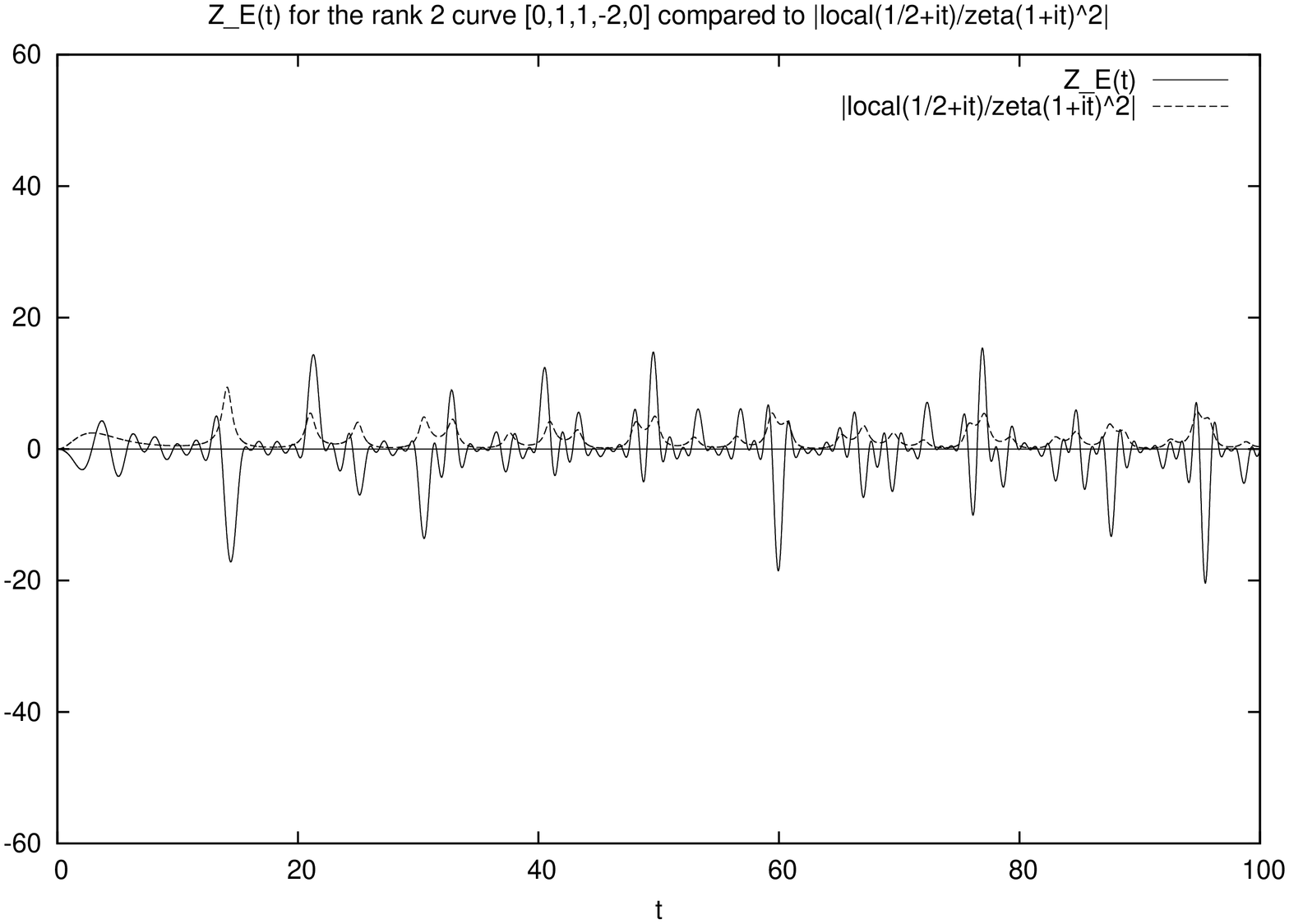,width=2.5in,angle=-90}
    \psfig{figure=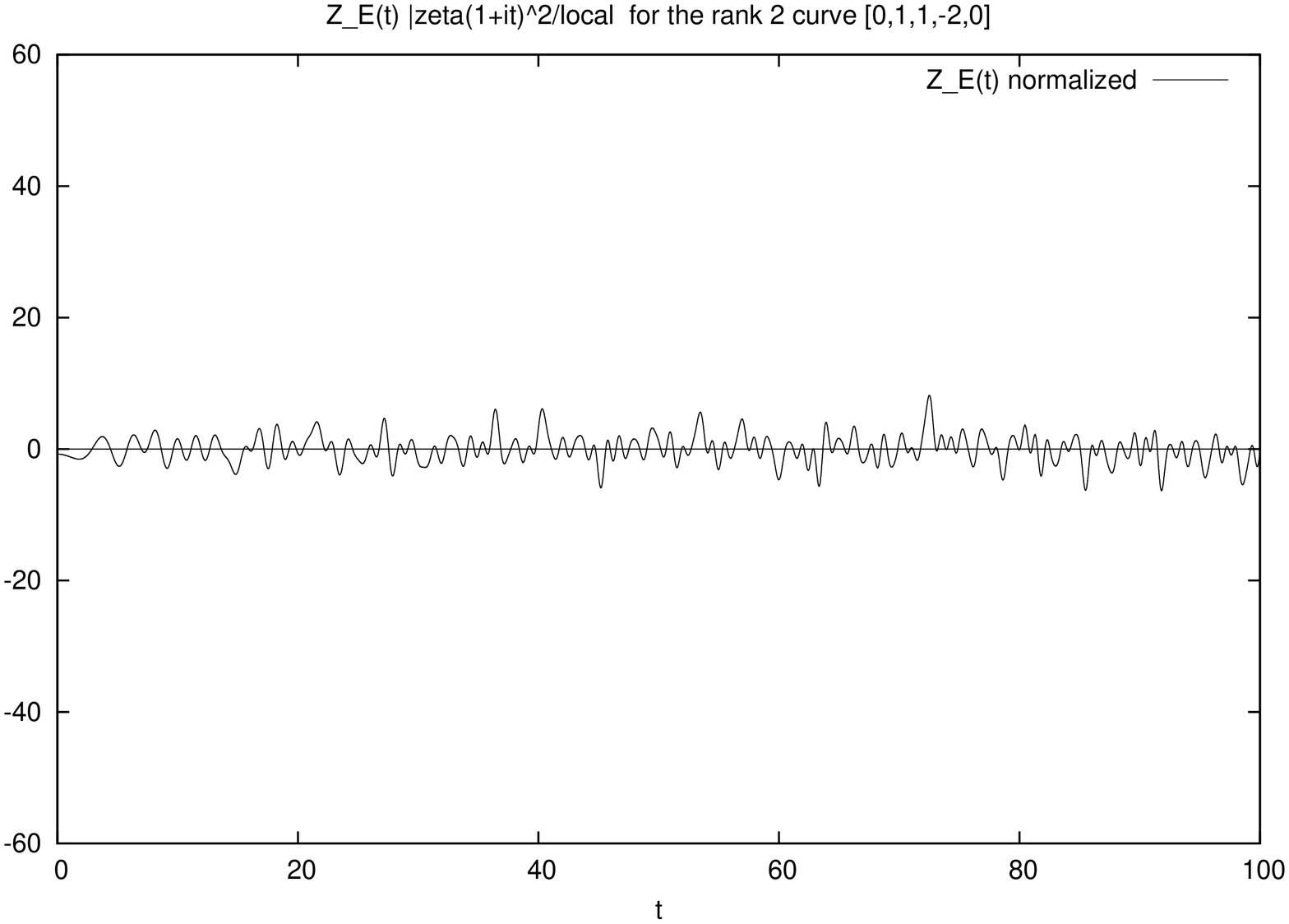,width=2.5in,angle=-90}
}
\centerline{
    \psfig{figure=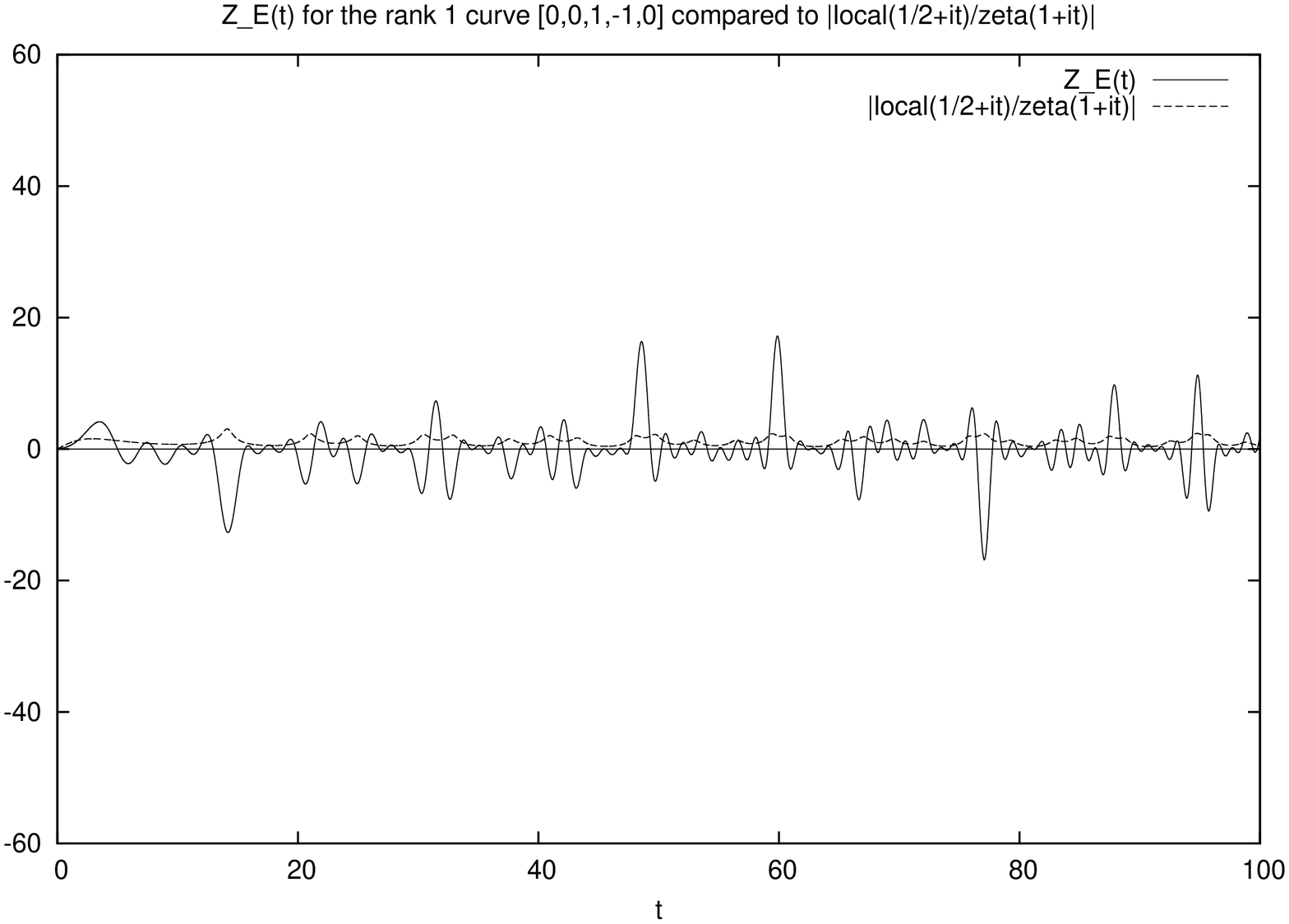,width=2.5in,angle=-90}
    \psfig{figure=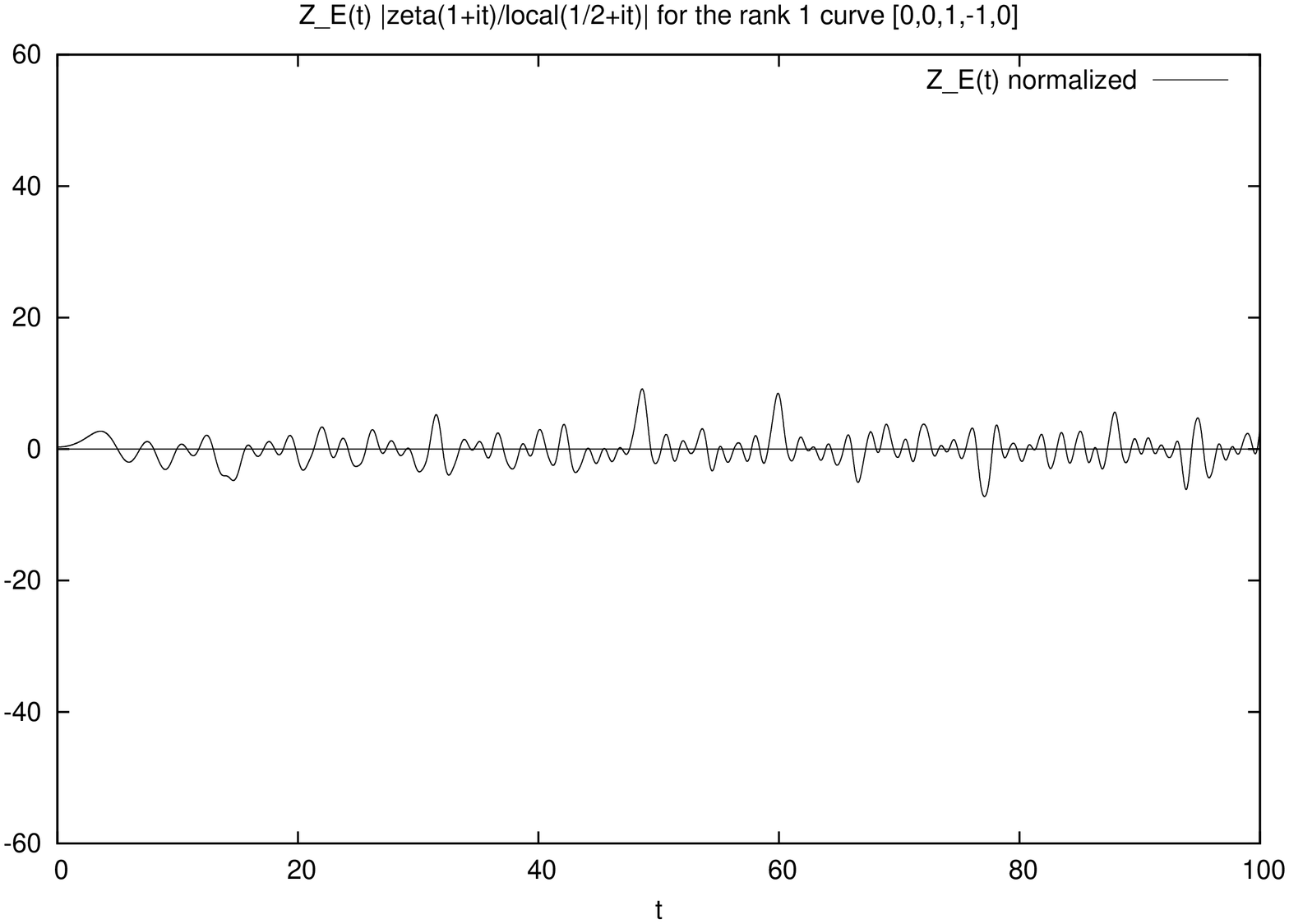,width=2.5in,angle=-90}
}
\caption
    {Same as previous page, but for $r=3,2,1$.}
\label{fig:graphs 3-1}
\end{figure}
}

Notice, in these plots, that $Z_E(t)$ tends to have a deficiency of zeros near
the $t=0$, as explained by the $r$ extra zeros that it acquires at $t=0$. This
feature should increase with $r$, but also dissipate, for given $r$, as the
conductor increases. See~\cite[Conjecture 1.1]{Mi} where this is conjectured for certain
families of elliptic curves, and also Section 4.2.2 of that paper which presents some
evidence in favour of this claim for $r=2$. See also~\cite{DHKMS} for a random matrix
model that explains this phenomenon.

Also notice the large gaps in Figures \ref{fig:spikes}, \ref{fig:graphs 7-5},
and \ref{fig:graphs 3-1} between zeros when $t$ is near the imaginary part of
the zeros of the zeta function, explained by the fact that $1/\zeta(1+it)^r$
tends to get large, especially initially, near these points.
We expect this phenomenon to also
dissipate as the conductor grows, and also as $t$ grows.

It would be interesting to see if these large gaps, near $t=0$ and, say, near
$t=14.134\ldots$, corresponding to the first zero of zeta, could be exploited
in analytic algorithms that make use, say, of the explicit formula.
See for example~\cite{BHK} for a novel algorithmic use of the explicit formula.

We also note, returning to~\eqref{eq:local factor expanded},
that the coefficient that accompanies the $p^{-2s}$ term,
\begin{equation}
    \notag
    a(p^2) = \frac{(a(p)^2-p)}{p} = \alpha^2 + \beta^2 + 1,
\end{equation}
is the Dirichlet coefficient for the prime $p$ of the symmetric square $L$-function
\begin{equation}
    \notag
    L_E(s,\text{symm}^2) = \prod_p (1-\alpha(p)^2 p^{-s})^{-1} (1-p^{-s})^{-1} (1-\beta(p)^2 p^{-s})^{-1}.
\end{equation}
This suggests that one should feel the presence of $L_E(1+2it,\text{symm}^2)$ when examining $L_E(1/2+it)$.
This is harder to see compared to the prominent high rank affect, but
$L_E(1+2it,\text{symm}^2)$ does show up in various averages of $L_E(1/2+it)$,
for example in lower terms of its moments or density of zeros, when averaged
over families of elliptic curves. See~\cite[Conj 2.1, Thm 2.2]{H} or~\cite{HMM}.

\subsection{Density of zeros for quadratic Dirichlet $L$-functions}

Interestingly, the Riemann zeta function on the one line appears in various
statistics of $L$-functions. One striking example concerns the density of zeros
of $L(s,\chi_d)$, where $\chi_d(n)=\br{\frac{d}{n}}$ is Kronecker's symbol.

\afterpage{
\begin{figure}[H]
    %XXXXXXXXXXX comment out for speed
    \centerline{
            \psfig{figure=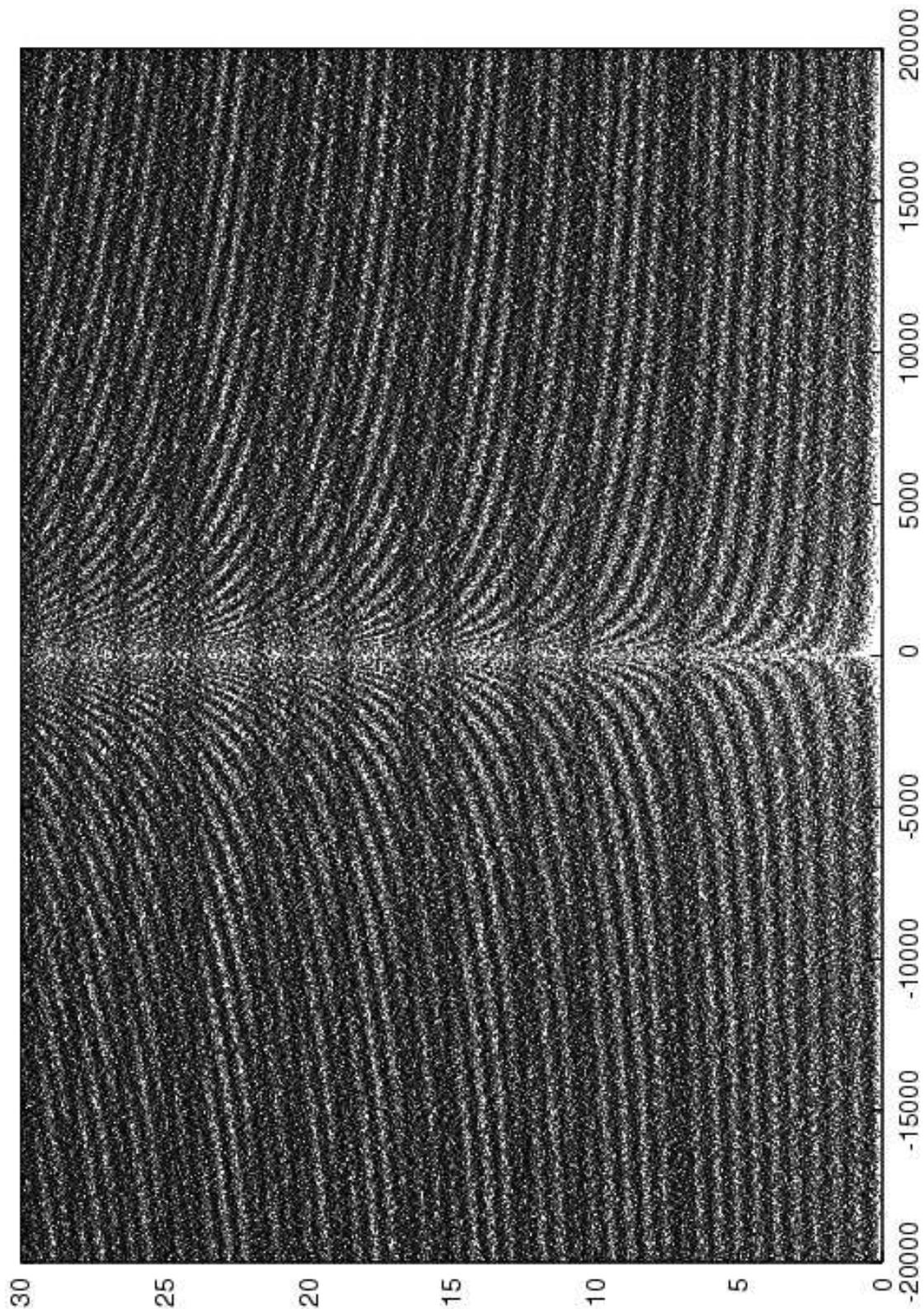,width=6in,angle=0}
    }
    \caption
    {Zeros of $L(s,\chi_d)$ with $\chi_d(n)=\left(\frac{d}{n} \right)$,
     the Kronecker symbol. We restrict $d$
     to fundamental discriminants $-20000< d < 20000$.
     The horizontal axis is $d$ and, for each
     $L(s,\chi_d)$, the imaginary parts of its zeros up to height $30$ are listed.}
    \label{fig:zeros real chi}
\end{figure}
}

Figure~\ref{fig:zeros real chi}, from~\cite{R3}, depicts the imaginary parts of the
non-trivial zeros of these
$L$-functions, for fundamental discriminants $d$, $|d|<20,000$.
We can observe the density of zeros fluctuating as one moves away from the real axis, and
also increasing slowly, as seen by the fact that the zeros tend to move towards
the real axis, as $|d|$ increases. The fact that the density increases
is easily explained by von Mangoldt's
formula for the number of zeros of $L(s,\chi_d)$ up to height $T$:
\begin{equation}
    \label{eq:von Mangoldt chi}
    \left|
    \{ \rho: L(\rho,\chi_d)=0, 0<\Re{\rho}<1, 0 \leq \Im{\rho} \leq T\}
    \right|
    \sim \frac{T\log(T|d|)}{2\pi}.
\end{equation}

Other features can be seen in the plot. First, the white band near the $x$-axis
indicates that the lowest zero of $L(s,\chi_d)$ repels away from the real axis.
We can also see the effect of secondary terms on this
repulsion. The lowest zero for $d>0$ tends to be higher than the lowest zero
for $d<0$. This turns out, as will be discussed below, to be related to the
fact that the $\Gamma$-factor in the functional equation for $L(s,\chi_d)$ is
$\Gamma(s/2)$ if $d>0$, but is $\Gamma((s+1)/2)$ when $d<0$.

Most relevant to our discussion are the slightly darker regions appearing in horizontal strips.
The first one occurs roughly at height $7.$, approximately half the imaginary part of
the first zero of $\zeta(s)$.
These horizontal strips are due to secondary terms in the density of zeros for this 
collection of $L$-functions
which include a term that is proportional to
$$
    \Re\frac{\zeta'(1+2it)}{\zeta(1+2it)}.
$$
See formula~\eqref{eq:density conrey snaith} below.
This tends to be large when $2t$ is near the imaginary part of a zero of the zeta function,
as can be seen from formula~\eqref{eq:zeta log diff} in the next section,
especially for smaller $t$, where the zeros are well spaced apart.

The fluctuating sand-dune like feature  is explained by the main term in the
density of zeros of $L(s,\chi_d)$.
Let $D(X)$ denote the set of fundamental discriminants up to $X$:
$$
   D(X) = \left\{
       \text{$d$ a fundamental discriminant : $|d| \leq X$}
   \right\},
$$
and let $f$ be smooth, rapidly decreasing, and having Fourier transform
supported in the interval $(-1,1)$.
\"Ozl\"uk  and Snyder proved~\cite{OS} that the average density of zeros of $L(s,\chi_d)$,
with test function $f$, satisfies
\begin{eqnarray}
    \label{eq:ozluk snyder}
    &&\lim_{X\rightarrow\infty} \frac{1}{|D(X)|} \sum_{d \in D(X)}
    \sum_{\gamma_d} f\left(\gamma_d \frac{\log |d|} {2\pi}
    \right) = \int_{-\infty}^{\infty} f(x) \bigg(1-\frac{\sin (2\pi
    x)}{2\pi x}\bigg) dx, \notag \\
\end{eqnarray}
where $\gamma_d$ runs over the imaginary parts of all the non-trivial zeros
of $L(s,\chi_d)$. The factor of $\log{|d|}/(2\pi)$ reflects the fact that the density
of non-trivial zeros in a fixed region near the real axis increases proportionately, in the $d$ aspect,
to this factor.
This scaling also has the effect of `zooming in' on the zeros close to the real
axis, in the sense that zeros satisfying $|\gamma_d| > \log(|d|)^{-1+\epsilon}$ contribute nothing to 
the limit, for any $\epsilon>0$.

\"Ozl\"uk  and Snyder also proved that the support condition can be relaxed to
the interval $(-2,2)$ if one assumes the GRH for $L(s,\chi_d)$.
Presumably, the theorem remains valid for
a wider class of test functions $f$, for example piecewise continuous
integrable functions $f: \R \to \R$.

\begin{figure}[htp]
    \centerline{
            \psfig{figure=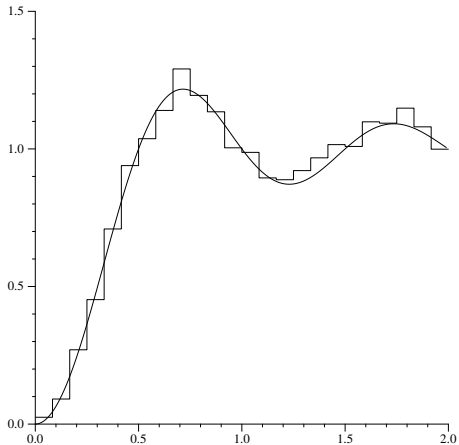,width=2.5in,angle=0}
    }
    \caption
    {Density of zeros of $L(s,\chi_d)$ for
     7243 prime values of $|d|$ lying in the interval $(10^{12},10^{12}+200000)$,
     compared against the  prediction, $1-\sin(2\pi x)/(2\pi x)$.}
    \label{fig:density zeros usp}
\end{figure}

In figure~\ref{fig:density zeros usp}, taken from~\cite{R}, we depict the
$1$-level density of the zeros of $L(s,\chi_d)$ for $7243$ prime $|d|$ lying in
the interval $(10^{12}, 10^{12}+200000)$.
Here the horizontal axis is divided into
small bins of width $1/10$, then count, on average, how many normalized zeros of
$L(s,\chi_d)$ lie in each bin, and find excellent agreement with the graph of $1-\sin(2\pi
x)/(2\pi x)$. See~\cite{R} for details regarding the normalization.

That the density is, to leading order,
$1-\sin(2\pi x)/(2\pi x)$ explains one of the basic features evident in
Figure~\ref{fig:zeros real chi}, specifically the fluctuating sand-dune like
regions in the plot. However, it fails to account for the more subtle features,
described above, that are evident in the figure. To do so requires the lower
terms in the density of zeros.

An approach has been developed by Conrey and Snaith~\cite{CS} for obtaining the
full asymptotic expansion of the density of zeros of $L(s,\chi_d)$. It is an
interesting application of precise conjectures of Conrey, Farmer, and
Zirnbauer, for moments of ratios of $L$-functions~\cite{CFZ}~\cite{CFS},
and is breathtaking in its detail.

To describe their formula, we let $g(z)$ be holomorphic throughout the strip
$|\Im z| <2$, real on the real line and even, and satisfy  $g(x)\ll 1/(1+x^2)$
as  $x\to \infty$. Subject to the `moments of ratios conjecture' and GRH for
$L(s,\chi_d)$, Conrey and Snaith proved:
{\small
\begin{eqnarray}
    \label{eq:density conrey snaith}
    &&\sum_{d\le X} \sum_{\gamma_d} g(\gamma_d)= \frac{1}{2\pi} \int_{-\infty}^\infty g(t)\sum_{d\le X}
    \bigg( \log \frac{d}{\pi} +\frac12 \frac {\Gamma'}{\Gamma}(1/4+it/2)\notag \\
    &&+\frac12 \frac {\Gamma'}{\Gamma}(1/4-it/2)+
    2\bigg(\frac{\zeta'(1+2it)}{\zeta(1+2it)} +A_D'(it;it)\notag \\
    &&-\left(\frac d \pi\right)^{-it}
    \frac{\Gamma(1/4-it/2)}{\Gamma(1/4+it/2)}
    \zeta(1-2it)  A_D(-it;it)\bigg)\bigg) ~dt
    +O(X^{1/2+\epsilon}),
\end{eqnarray}
}
where
\begin{eqnarray}
    A_D(-r;r)=
    \prod_p\left(1-\frac{1}{(p+1)p^{1-2r}}-\frac{1}{p+1}\right)\left(1-\frac{1}{p}\right)^{-1},
\end{eqnarray}
and
\begin{eqnarray}
    A'_D(r;r)=\sum_p \frac{\log p}{(p+1)(p^{1+2r}-1)}.
\end{eqnarray}
Again, this formula presumably continues to hold for piecewise continuous integrable
functions $g:\R \to \R$.

\afterpage{
\begin{figure}[H]
    \centerline{
            \psfig{figure=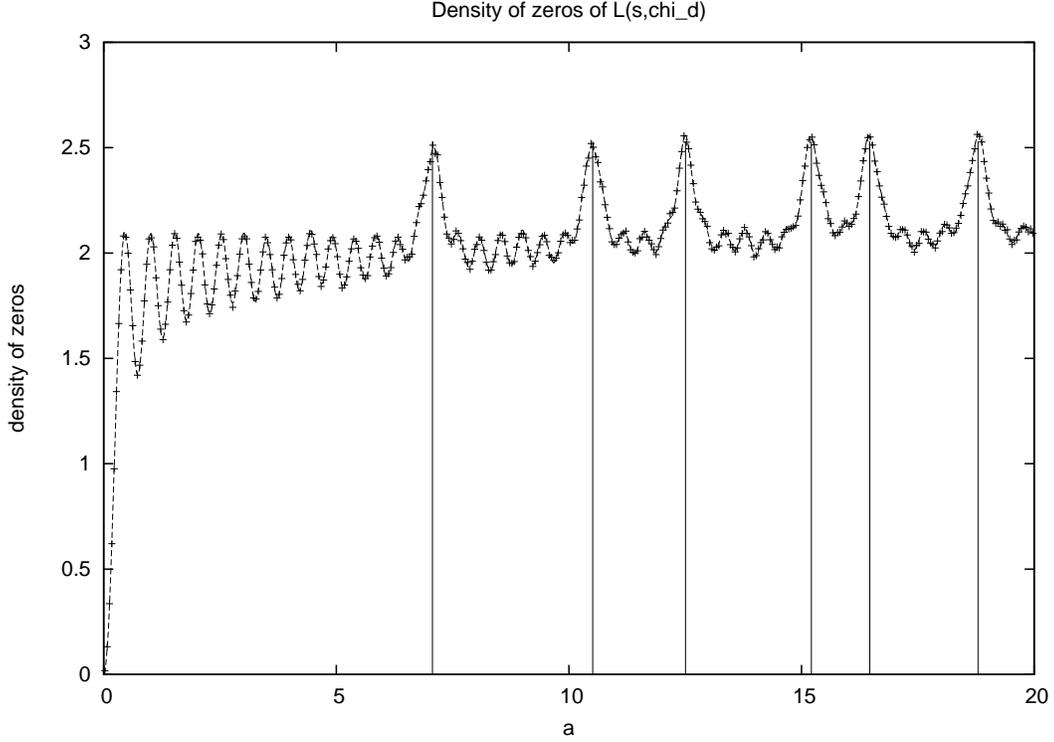,width=4in,angle=270}
    }
    \caption{The density of zeros of $L(s,\chi_d)$ for fundamental discriminants
    $0<d<10^6$ (histogram), compared to the prediction of Conrey and Snaith
    (dashed curve). The vertical lines mark half the imaginary parts of the
    zeros of zeta.}
    \label{fig:one level cs}
\end{figure}
}

We compare both sides of~\eqref{eq:density conrey snaith} in
Figure~\ref{fig:one level cs}, with $X=10^6$. For each $0 < d < 10^6$, we
computed the first 100 zeros above the real axis of $L(s,\chi_d)$ using the author's
{\tt lcalc} package~\cite{R2}. We then let
$h(x) = \chi_{[a,a+1/20)}(x)$, for $a=0, .05, .1, .15, \ldots, 19.95$, i.e.
indicator functions of intervals of width $1/20$, and take $g(x)
=(20/303957)\times (h(x)+h(-x))/2$, so that $g(x)$ is even. The factor of $20$
is to account for the width of the interval being $1/20$, while the $303957$ is
the number of positive fundamental discriminants $0 < d < 10^6$, so that we are averaging.
For each value of $0\leq a \leq 19.95$, we display the lhs of~\eqref{eq:density
conrey snaith} as `plus' marks, with each `+' being centred on the interval
$[a,a+1/20)$, i.e. having horizontal coordinate $a+1/40$. For the rhs, we plot, as a dashed curve,
the following approximation:
{\small
\begin{eqnarray}
    &&\frac{1}{2\pi} \frac{1}{303957} \sum_{d\le 10^6} \Re
    \bigg( \log \frac{d}{\pi} + \frac {\Gamma'}{\Gamma}(1/4+ia/2) + \notag \\
    %+\frac12 \frac {\Gamma'}{\Gamma}(1/4-ia/2)+ \notag \\
    &&2\bigg(\frac{\zeta'(1+2ia)}{\zeta(1+2ia)} +A_D'(ia;ia)
    -\left(\frac d \pi\right)^{-ia}
    \frac{\Gamma(1/4-ia/2)}{\Gamma(1/4+ia/2)}
    \zeta(1-2ia)  A_D(-ia;ia)\bigg)\bigg),\notag \\
\end{eqnarray}
}
i.e. the real part of the integrand without the factor of $20$.
The agreement is beautiful.
%Because the width of our bins is small, $1/20$,
%we plot, for the rhs the values of the integrand, including the factor of
%$1/(2\pi)$ in front of the integral, as a dashed curve. Because the width of
%our bins is small, $1/20$, the integrand does a good job

Note that the formula is stated for $d>0$. A nearly identical formula holds for
$d<0$, but with the $1/4+it/2$ being replaced by $3/4+it/2$ throughout the
formula, reflecting the $\Gamma( (s+1)/2)$ factor that appears in the functional
equation of $L(s,\chi_d)$ when $d<0$. Note that
$\Gamma'(1/4)/\Gamma(1/4) = -4.228\ldots$, whereas
$\Gamma'(3/4)/\Gamma(3/4) = -1.085\ldots$, which explains, for example, the fact
that the first zero of $L(s,\chi_d)$ tends to lie slightly closer to the real
axis for $d<0$ as compared to $d>0$.

One can recover the formula of \"Ozl\"uk  and Snyder~\eqref{eq:ozluk snyder} by
applying~\eqref{eq:density conrey snaith} to the function
$g(t)=f(t\log(X)/(2\pi))$, substituting $\tau=t\log X/(2\pi)$ in the integral
on the rhs of~\eqref{eq:density conrey snaith}, identifying the leading order
terms, and also using the fact that $\log|d| \sim \log{X}$ for most $|d|\leq
X$. See~\cite{CS} for the details.

For completeness, we mention that the function $1-\sin(2\pi x)/(2\pi x)$ is the
limiting density function for the eigenangles, suitably normalized, of large
unitary symplectic matrices, consistent with predictions of Katz and
Sarnak~\cite{KS}~\cite{KS2}. This agreement has been found to persist, more
generally, for a higher dimensional analogue involving the density of
$r$-tuples of zeros~\cite{R}~\cite{ERR}.

\section{$\zeta(1+it)$ in the pair correlation of zeta zeros}

Another interesting connection with the Riemann zeta function on the one line
occurs in the pair correlation of the non-trivial zeros of the zeta function.
It was here that the subtle influence that the Riemann zeta function on the one
line asserts on $L$-functions was first discovered by Bogomolny and Keating
~\cite{BK}~\cite{BK2}~\cite{BBLM}.

The Riemann Hypothesis states that the non-trivial zeros of the Riemann zeta
function have real part equal to one half. We thus write a typical non-trivial
zero of $\zeta$ as $1/2 + i\gamma$,
and assume that the $\gamma$'s are real. The zeros come in conjugate pairs, so we
label those above the real axis as: $0 < \gamma_1 \leq \gamma_2 \leq \ldots$.

Montgomery, was the first to study the vertical distribution of the zeros~\cite{Mo}.
He considered the pair correlation of the zeros, a basic statistic that
measures how much the zeros know about one another.

Let
$
    N(T)
$
denote the number of non-trivial zeros of $\zeta(s)$ with $0< \Im(s)\leq T$.
A theorem of von Mangoldt states that
\begin{equation}
    N(T) \sim \frac{T \log{T}}{2\pi},
    \label{eq:von Mangoldt}
\end{equation}
as $T \to \infty$. Thus, by scaling each $0<\gamma\leq T$ by  $\log(T)/2\pi$,
they are spaced apart by one, on average.

Let $\alpha < \beta$. Montgomery conjectured that
\begin{eqnarray}
    \frac{1}{N(T)}
    \left| \left\{ 1 \leq i \neq j \leq N(T): (\gamma_j - \gamma_i) \tfrac{\log{T}}{2\pi}
        \in [\alpha,\beta ] \right\} \right| \notag \\
    \sim
    \int_\alpha^\beta
    \left( 1 - \left(\frac{\sin \pi t}{\pi t} \right)^2 \right) dt,
    \label{eq:montgomery conj}
\end{eqnarray}
as $T \to \infty$. Notice that the integrand is small when $t$ is near 0, so that
zeros of $\zeta$ tend to repel away from one another.
Odlyzko carried out extensive numerical tests
of Montgomery's conjecture~\cite{O}~\cite{O2}~\cite{O3}.

The factor of $\log(T)/2\pi$ in~\eqref{eq:montgomery conj} is to account for
the fact that the zeros become more dense as one moves away from the real axis
in accordance with~\eqref{eq:von Mangoldt}.

%\begin{frame}
Using the explicit formula that connects zeros of zeta to prime powers,
Montgomery was able to prove that
\begin{eqnarray}
    \label{eq:montgomery theorem}
    \frac{1}{N(T)}
    \left| \left\{ 1 \leq i \neq j \leq N(T): f\left((\gamma_j - \gamma_i) \tfrac{\log{T}}{2\pi} \right)
    \right\} \right|
    \to
    \int_{-\infty}^\infty f(t) \left( 1 - \left(\frac{\sin \pi t}{\pi t} \right)^2 \right) dt, \notag \\
\end{eqnarray}
as $T \to \infty$,
for smooth and rapidly decaying functions $f$ satisfying the stringent restriction that $\hat{f}$
be supported in $(-1,1)$.
Rudnick and Sarnak~\cite{RS} generalized this to any primitive $L$-function (assuming
a weak form of the Ramanujan conjectures in the case of higher degree $L$-functions).
They also gave a smoothed version of the above theorem in the case that RH is false.

%\end{frame}

%%-------------------begin frame: Dyson -----------------------
%%\begin{frame}
%        On a visit by Montgomery to the Institute for Advanced Study,
%        Freeman Dyson pointed out that large unitary matrices have the same
%        pair correlation. Let
%        $$
%            e^{i \theta_1}, e^{i \theta_2}, \ldots, e^{i \theta_N}
%        $$
%        be the eigenvalues of a matrix in $\text{U}(N)$, sorted so that
%        $$
%            0 \leq \theta_1 \leq \theta_2 \ldots \leq \theta_N < 2 \pi.
%        $$
%        Normalize the eigenangles
%        $$
%             \tilde{\theta}_i = \theta_i N / (2\pi)
%        $$
%        so that $\tilde{\theta}_{i+1}-\tilde{\theta}_i$ equals one on average.
%%\end{frame}
%%\begin{frame}

This result connects, statistically, the zeros of $\zeta(s)$, to eigenvalues
of large unitary matrices. If $\exp(i \theta_j)$ are the eigenvalues of
a matrix in $U(N)$, $j=1,\ldots, N$, and
$$
    0 \leq \theta_j < 2 \pi,
$$
then, a classic result in random matrix theory~\cite{M} asserts that
$$
    \frac{1}{N}
    \left| \left\{  1 \leq i \neq j \leq N:
    (\theta_j-\theta_i)\frac{N}{2\pi} \in [\alpha,\beta ] \right\} \right|
$$
equals, when averaged according to Haar measure over $\text{U}(N)$ and letting
$N \to \infty$,
$$
    \int_\alpha^\beta
    \left( 1 - \left(\frac{\sin \pi t}{\pi t} \right)^2 \right) dt.
$$
%\end{frame}
%%------------------ end frame -------------------------------------

Interestingly, very precise lower terms have been conjectured for the pair correlation.
These were first described by Bogomolny and Keating. They used the Hardy-Littlewood conjecture
for the asymptotic number of prime pairs with given difference to estimate `off-diagonal'
contributions~\cite{BK}~\cite{BK2}.

We detail the conjectured lower terms in the pair correlation as described by
Conrey and Snaith~\cite{CS}. As in the previous section, Let $g(z)$ be
holomorphic throughout the strip $|\Im z| <2$, real on the real line and even,
and satisfy  $g(x)\ll 1/(1+x^2)$ as  $x\to \infty$. Subject to the `moments of
ratios conjecture' and Riemann Hypothesis for the zeta function, Conrey and Snaith proved:

\begin{eqnarray}
    \notag
    &&\sum_{ 1 \leq i \neq j \leq N(T)}  g(\gamma_j-\gamma_i)
     =\frac{1}{(2\pi)^2}\int_0^T
    \int_{-T}^T g(r) \bigg( \log^2 \frac{t}{2\pi} +2\bigg(\left(\frac{\zeta'}{\zeta}\right)'(1+ir)\\
     && \qquad \qquad
    + \left(\frac{t}{2\pi}\right)^ {-ir}
    \zeta(1-ir)\zeta(1+ir)A(ir)-B(ir)\bigg)\bigg) ~dr ~dt
    +O(T^{1/2+\epsilon}), \notag \\
    \label{eq:conrey snaith lower terms}
\end{eqnarray}
where the integral over $r$ is regarded as principal valued near $r=0$,
\begin{eqnarray}
A(\eta)=\prod_p\left(1-\frac{1}{p^{1+\eta}}\right)
\left(1-\frac 2 p +\frac{1}{p^{1+\eta}}\right)\left(1-\frac 1 p \right)^{-2},
\end{eqnarray}
and
\begin{eqnarray}
    B(\eta)=\sum_p \left(\frac{\log p}{(p^{1+\eta}-1)}\right)^2.
\end{eqnarray}
Presumably, formula~\eqref{eq:conrey snaith lower terms} continues to hold,
assuming RH, for a wider class of test functions $f$, for example piecewise continuous
integrable functions $g: \R \to \R$.

Conrey and Snaith also showed that one recovers~\eqref{eq:montgomery theorem}
by letting  $g(x)=f(x \frac{\log{T}}{2\pi})$, and
substituting $y=r\frac{\log{T}}{2\pi}$ in the inner integral above.
%\begin{eqnarray}
%    \sum_{0 < \gamma \neq \gamma'\le T} f((\gamma-\gamma')\tfrac{\log{T}}{2\pi})
%    \sim
%    N(T)
%    \int_{-\infty}^{\infty}
%    f(y)
%    \left( 1 - \left(\frac{\sin \pi y}{\pi y} \right)^2 \right). \notag
%\end{eqnarray}

%Recently, Ford and Zaharescu~\cite{FZ} have obtained, unconditionally, some of the lower terms in
%the above conjecture.

Notice the term $(\zeta'/\zeta)'(1+ir)$ which tends to be large, in magnitude, at least initially,
when $r$ is close to the imaginary part of a non-trivial zero of $\zeta$. This can be
explained by taking the derivative of the well known formula~\cite[Chapter 12]{D}
\begin{equation}
    \label{eq:zeta log diff}
    \frac{\zeta'(s)}{\zeta(s)} = C - (s-1)^{-1} - \Psi(s/2+1)/2
    + \sum_\rho \left( (s-\rho)^{-1} + \rho^{-1} \right)
\end{equation}
where $C$ is a constant, the sum is over the nontrivial zeros $\rho = 1/2+i\gamma$
of $\zeta$, and $\Psi(z)=\Gamma'(z)/\Gamma(z)$. On differentiating, the sum over $\rho$ becomes
\begin{equation}
    -\sum_\rho (s-\rho)^{-2}.
\end{equation}
For a given $\rho=1/2+i\gamma$, and $s=1+ir$, the corresponding term in the above sum is
largest in magnitude when $r=\gamma$, in which case $(s-\rho)^{-2}=4$. Furthermore, the first
few non-trivial zeros of $\zeta$ are well spaced apart, so that, for smaller $r$, the influence
of these first few zeros is felt quite distinctly. On differentiating,
the other terms in~\eqref{eq:zeta log diff} contribute little as $r$ grows, as can be seen
using the asymptotic formula $\Psi'(z) \sim 1/z$.

This is illustrated in
Figure~\ref{fig:3 plots one line} which depicts $|(\zeta'/\zeta)'(1+ir)|$.
Notice the peaks of height roughly 4 near the first few $\gamma$, which are indicated
by vertical lines.
We also plot, because of its relevance to the sections of this paper,
the related graphs of $|(\zeta'/\zeta)(1+ir)|$ and of $|1/\zeta(1+ir)|$.

\afterpage{
\begin{figure}[H]
    \centerline{
        \psfig{figure=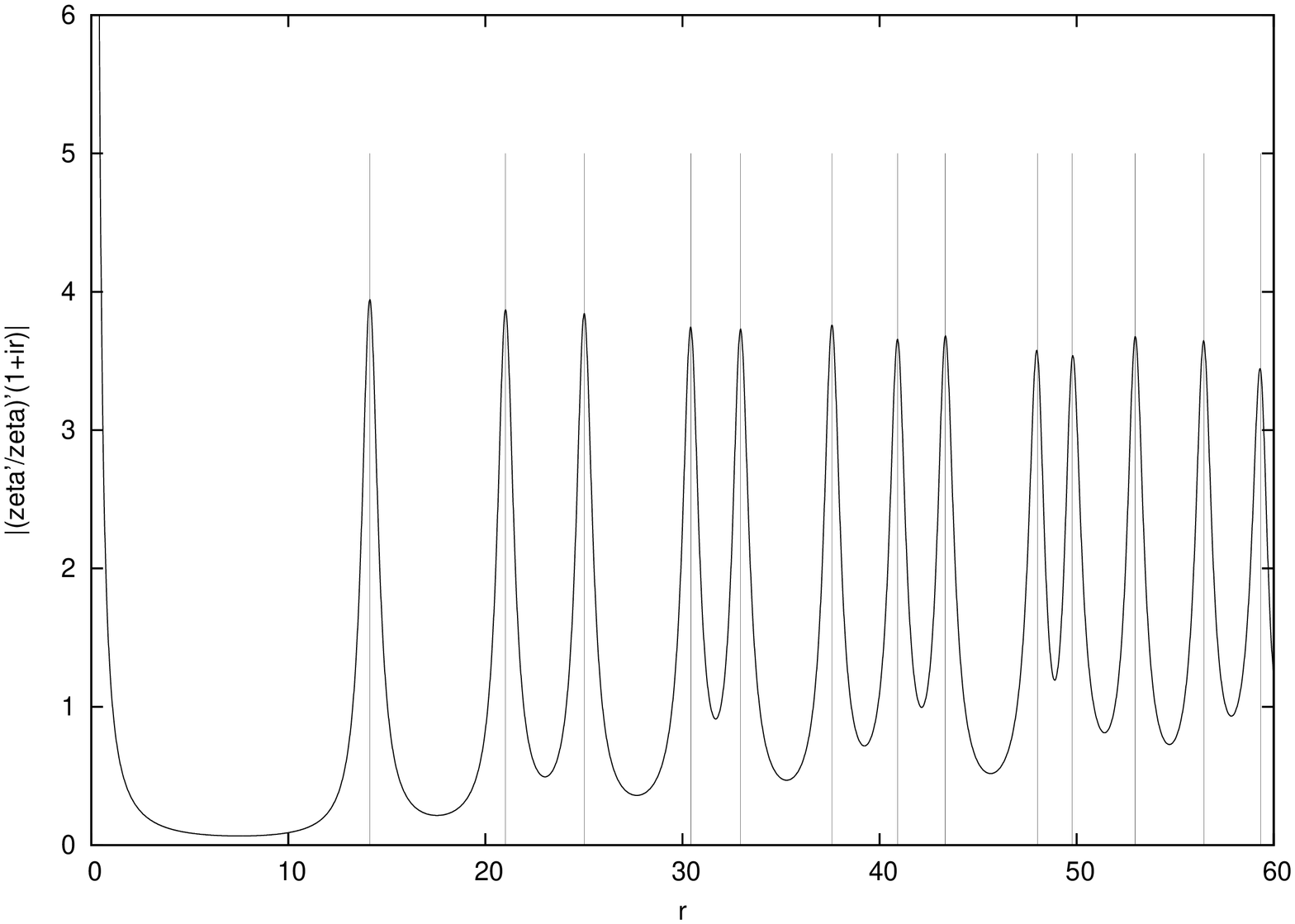,width=2.8in,angle=270}
    }
    \centerline{
        \psfig{figure=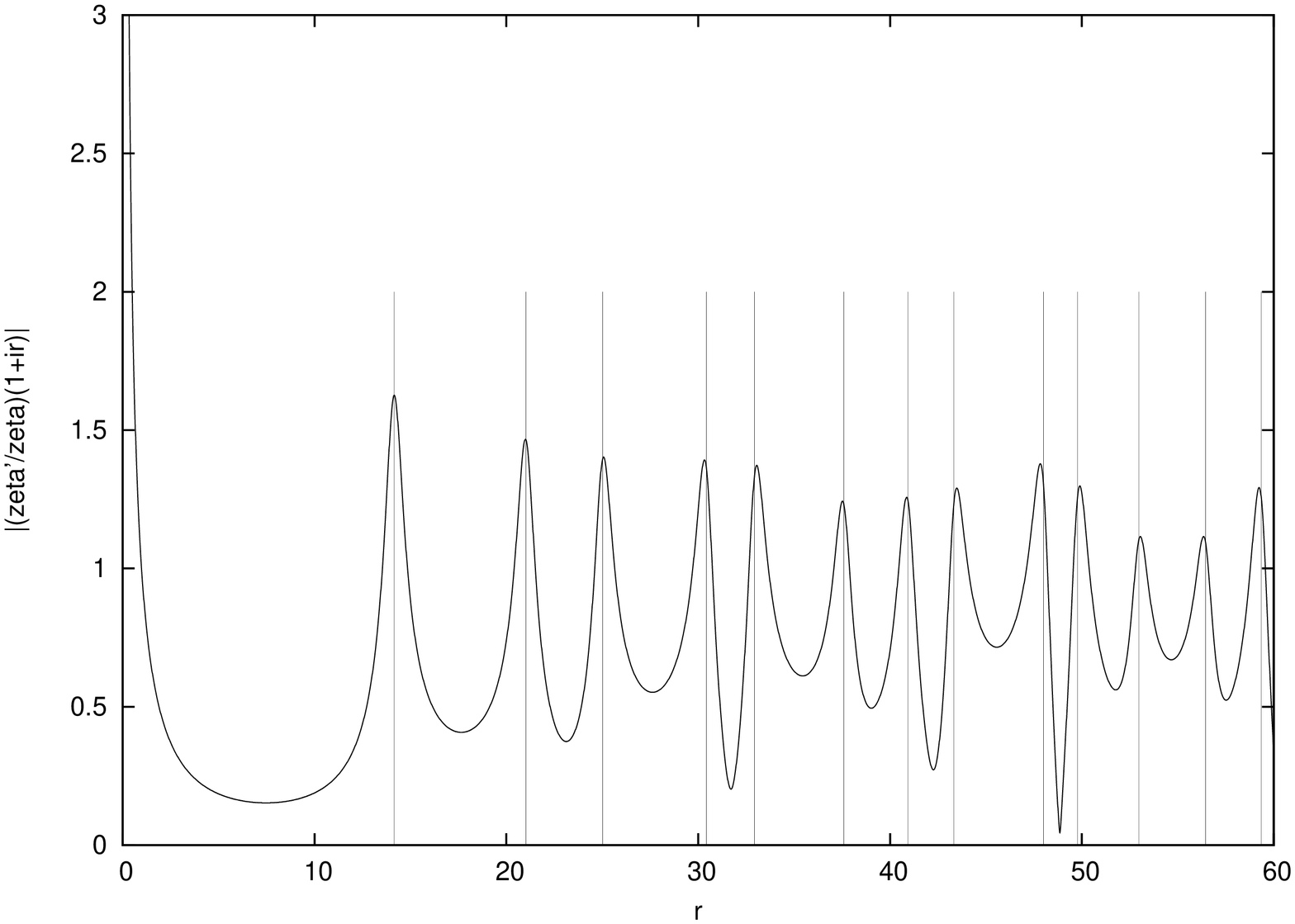,width=2.8in,angle=270}
    }
    \centerline{
        \psfig{figure=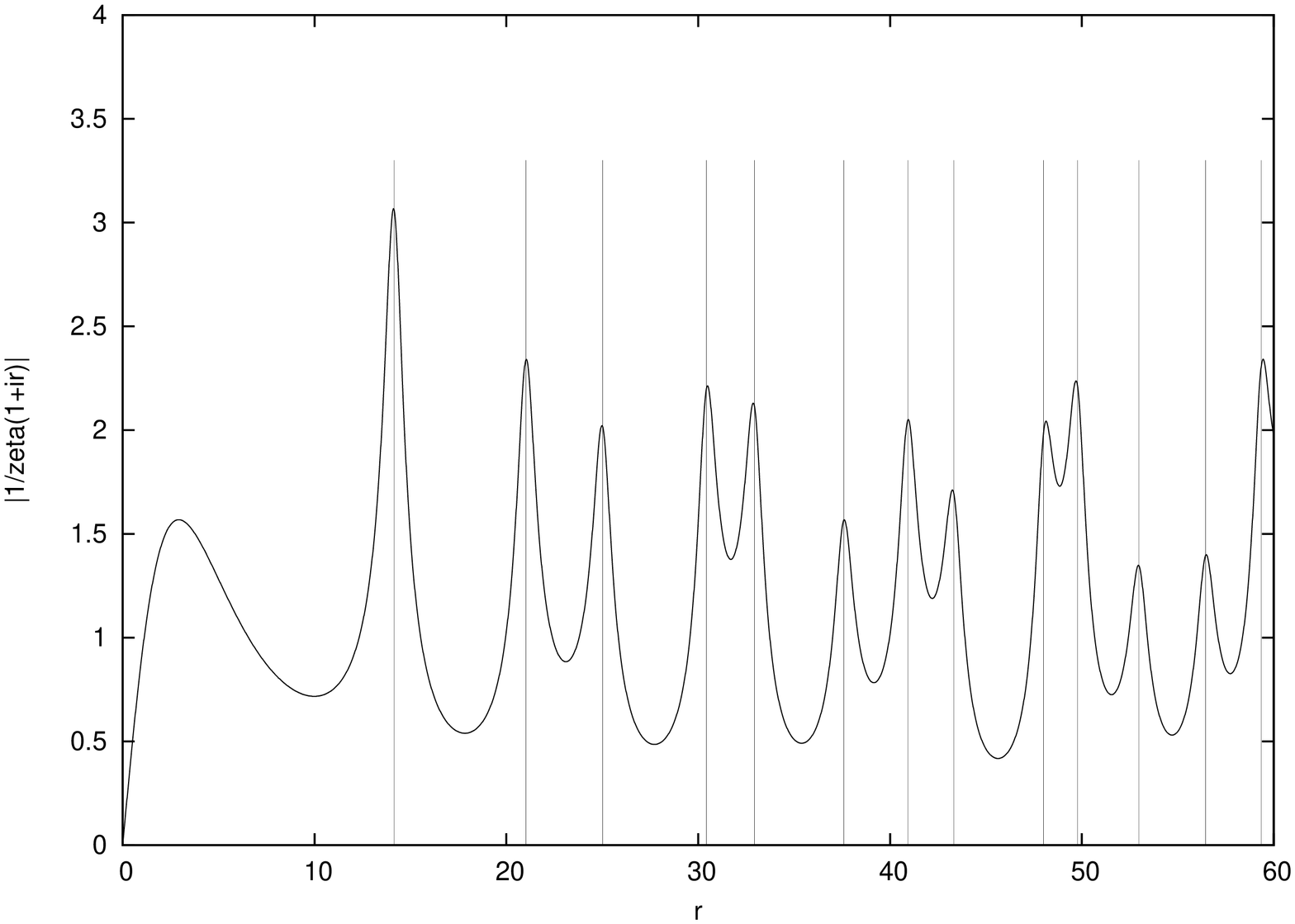,width=2.8in,angle=270}
    }
    \caption{Plots of $|(\zeta'/\zeta)'(1+ir)|,|(\zeta'/\zeta)(1+ir)|, |1/\zeta(1+ir)|$. The vertical
    lines mark the imaginary parts of the first few non-trivial zeros of $\zeta$.}
    \label{fig:3 plots one line}
\end{figure}
}

%Odlyzko carried out extensive computations
%Odlyzko data: $2\times 10^8$ zeros of zeta near the $10^{23}$rd zero.
%Pair correlation from data, bins of size $.01$, versus $1 - \sin(\pi t)^2/(\pi t)^2$.
%
%\begin{figure}[H]
%    \centerline{
%        \psfig{figure=odlyzko_pair_correlation.ps,width=2.7in,angle=270}
%        \psfig{figure=odlyzko_pair_correlation2.ps,width=2.7in,angle=270}
%    }
%    \caption
%    {Odlyzko's comparison of the pair correlation picture for $2\times 10^8$ zeros of $\zeta(s)$ near the
%     $10^{23}$rd zero to Montgomery's prediction. The second graph shows the difference between
%     the histogram in the first graph and $1 - \left((\sin \pi t)/(\pi t) \right)^2$.
%     In the interval displayed, the two agree to within about $.002$.
%    }
%    \label{fig:odlyzko1}
%\end{figure}

We end with a plot, reprinted from Snaith's paper~\cite{Sn}, that compares both sides
of~\eqref{eq:conrey snaith lower terms} for the first $100,000$ non-trivial
zeros of the zeta function, and, for $g$, many small bins of width $1/40$.
\begin{figure}[H]
    \centerline{
        \psfig{figure=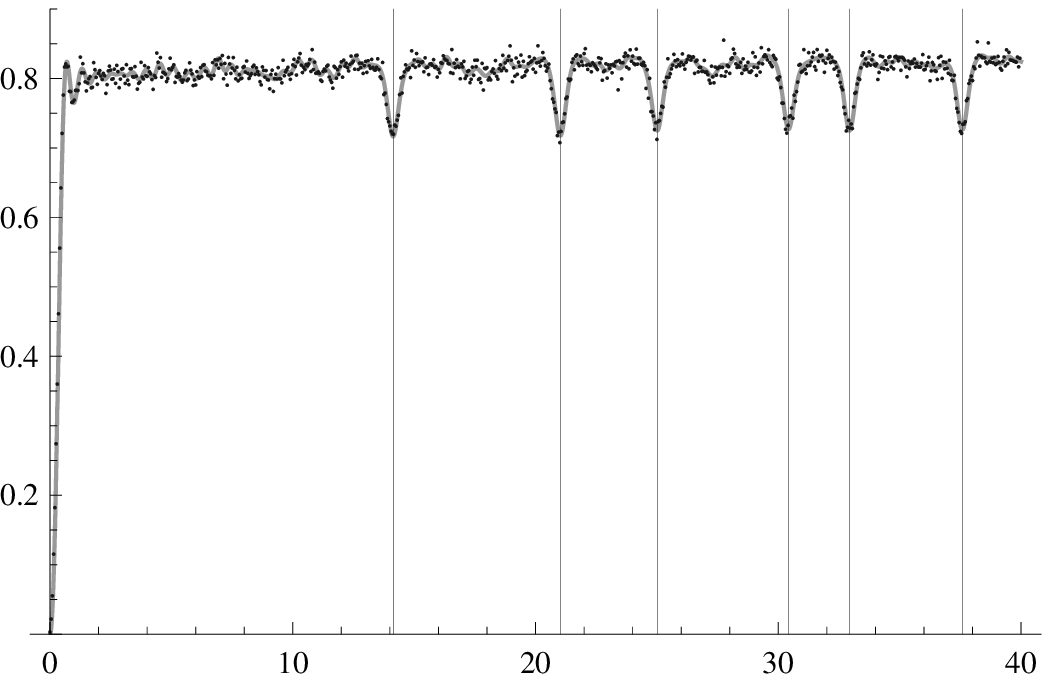,width=4in,angle=0}
    }
    \caption{A comparison of the pair correlation of the first $100,000$ zeros of $\zeta(s)$
    and the prediction given in~\eqref{eq:conrey snaith lower terms}. The vertical lines
    mark the imaginary parts of the first few non-trivial zeros of the zeta function.
    Courtesy of Nina Snaith~\cite{Sn}.}
\end{figure}
%------------------ end frame -------------------------------------

\subsection*{Acknowledgment}

This work was supported by the NSF Focused Research Group grant DMS 0757627
and by the author's NSERC Discovery grant. This paper grew out of a letter that
the author wrote in 2008. I would like to thank John Cremona, Noam Elkies,
David Farmer, Andrew Granville, Steven J. Miller, Mark Watkins, and the referee for feedback.

% ------------------------------------------------------------------------
\bibliographystyle{amsplain}

% ------------------------------------------------------------------------

\end{document}